\documentclass[12pt]{article}
\usepackage{amssymb}
\usepackage{amsmath}
\usepackage{amsopn}
\numberwithin{equation}{section}

\newtheorem{theorem}{Theorem}
\newtheorem{proposition}{Proposition}[section] 
\newtheorem{corollary}[proposition]{Corollary}
\newtheorem{lemma}[proposition]{Lemma}

\let\ds=\displaystyle

\def\coun{\varepsilon}
\def\bn{\begin{equation}}
\def\ed{\end{equation}}

\def\<{\langle}
\def\>{\rangle}
\newcommand{\CC}{{\mathbb C}}

\newcommand{\ZZ}{{\mathbb Z}}

\def\r#1{(\ref{#1})}
\let\rf=\r
\def\ot{\otimes}
 \def\d{\delta}

\def\sk#1{\left(#1\right)}
\def\U{\overline{U}}
\def\i{\iota}

\def\Pfp{{P}^+}
\def\Pfm{{P}^-}

\def\tPfp{{\tilde P}^+}

\def\Uqgln{U_q(\widehat{\mathfrak{gl}}_N)}
\def\Uqgl#1{U_q(\widehat{\mathfrak{gl}}_{#1})}

\def\Uqsl2{U_q(\widehat{\mathfrak{sl}}_2)}

\def\ff{{F}}
\def\EE{{\rm E}}\def\FF{{\rm F}}

\def\w{\mathbf{w}}
\def\prodl{\mathop{\overleftarrow\prod}\limits}

\def\ct{\mathbb{T}}
\def\End{\textrm{End}}

\def\bbb{\mathbb{B}}

\def\ccR{\mathbb{R}}
\def\si{\sigma}
\def\RR{{\rm R}}
\def\LL{{\rm L}}
\def\ord{\prec}
\def\E{{\rm E}}
\def\F{{\mathcal{F}}}
\def\tF{\tilde{\mathcal{F}}}

\def\tSym{\overline{\rm Sym}^{\ (q)}}
\def\qSym{\overline{\rm Sym}}
\def\ele{{\sf E}}
\def\tele{\tilde{\sf E}}
\def\hba{\nu}
\def\Uqqq{U_q(\mathfrak{gl}_N)}
\def\weig{\Lambda}

\def\lll{l}
\def\rr{r}
\def\ss{s}

\def\admis#1{[[\bar #1]]}
\def\seg#1#2{[#2,#1]}
\def\segg#1{[#1]}
\def\pp{p}
\def\ppp{\tilde{p}}
\def\Ff{f}
\def\acc{G}
\let\qsym=\beta
\def\tLL{{\tilde{\rm L}}}

\def\tX{{\tilde X}}
\def\tY{{\tilde Y}}
\def\tcZ{{\tilde{\mathcal{Z}}}}
\def\tcW{{\tilde{\mathcal{W}}}}
\def\tbfw{{\tilde{\bf w}}}
\def\coA{A^\perp}

\begin{document}

\begin{center}

\hfill ITEP-TH-66/06\\
\bigskip
{\Large\bf A computation of  an universal weight function for
the quantum affine algebra $\Uqgln$}
\par\bigskip\medskip
{\bf
Sergey Khoroshkin$^{\star}$\footnote{E-mail: khor@itep.ru}\ \  and \  \
Stanislav Pakuliak$^{\star\bullet}$\footnote{E-mail: pakuliak@theor.jinr.ru}}\ \
\par\bigskip\medskip
$^\star${\it Institute of Theoretical \& Experimental Physics, 117259
Moscow, Russia}
\par\smallskip
$^\bullet${\it Laboratory of Theoretical Physics, JINR,
141980 Dubna, Moscow reg., Russia}\\
\bigskip

\end{center}

\thispagestyle{empty}

\begin{abstract}
We compute an universal  weight function (off-shell Bethe
vectors) in any representation with a weight singular vector
of  the quantum affine algebra $U_q(\widehat{\mathfrak{gl}}_N)$
 applying the method of projections of Drinfeld  currents developed in \cite{EKhP}.
\end{abstract}


\section{Introduction} ${}$

 In  \cite{EKhP} a new method for the construction of  weight functions,
also known as off-shell Bethe vectors, was suggested. In this paper
we apply the constructions of \cite{EKhP} and \cite{KPT} to the
quantum affine algebra $\Uqgln$ and compute explicitly off-shell
Bethe vectors in representations with a weight singular vector.

Off-shell Bethe vectors in integrable models associated with the Lie algebra
$\mathfrak{gl}_N$ have appeared in \cite{KR83} in the framework of the algebraic
nested Bethe ansatz.
Except for $N=2$, off-shell Bethe vectors are defined  inductively.
They are functions of several complex variables. If the variables  satisfy the Bethe
ansatz equations, the Bethe vectors are eigenvectors of the transfer matrix of
the system. The construction of \cite{KR83} was developed in \cite{VT}, where weight
functions were defined as certain matrix elements of monodromy operators. Recently, the
construction of \cite{VT} was used in \cite{TV3} for the calculation of nested Bethe
 vectors in evaluation modules and their tensor products.

The approach of \cite{EKhP}
 is based on the ``new realization'' of quantum affine
algebras \cite{D88} and its connection to fundamental coalgebraic properties of
 weight functions,
found in \cite{VT}. The key role is played by certain
projections to the intersection of Borel subalgebras of
the quantum affine algebra of different kind, introduced in \cite{ER}.
 It is shown in
\cite{KP,EKhP} that acting by a projection of a product of Drinfeld currents
on highest weight vectors of irreducible finite-dimensional representations
of $\Uqgln$ one obtains a collection of rational functions with the required
comultiplication properties, that is, a weight function.

In this paper we perform a calculation of those projections and present
 explicit expressions for weight functions of the quantum affine algebra
  $\Uqgln$ in modules with a weight singular vector.
We present the result in two different ways. First, we write the
weight function as a polynomial  over Gauss coordinates of defining
$L$-operators with rational dependence of parameters, see
 eq \rf{Wt2}; this is equivalent to certain Cauchy type
integral over Drinfeld currents, see Theorem 1; eq. \rf{cc-pr1} and
eq. \rf{pr-st-inv}. Second, we express the weight function  as a
polynomial with rational functional coefficients over matrix
elements of the same $L$-operators (Theorem 2). As an application we
compute off-shell Bethe vectors in evaluation modules.

We compare our results with calculation of \cite{TV3} and verify
  a variant of the conjecture of \cite{KPT} about the coincidence of
  the two  constructions of weight functions. For such a verification we
  need to adjust our construction to
the data of \cite{TV3}, were the different $R$-matrix is used.
Two types of the recurrence relations for the universal off-shell Bethe vectors
related to two different embedding of $\Uqgl{N-1}$  into $\Uqgln$
was obtained in \cite{TV3}. In the forthcoming paper \cite{OPS} we will explain
these two types of the relations for the universal weight function using
two isomorphic current realizations of $\Uqgln$.

The paper is organized as follows. In section 2 we recall  the two
descriptions of $\Uqgln$: by means of the fundamental $L$-operators and
the so called 'new realization' of Drinfeld. In  Section 3 we describe, following
 \cite{KPT}, the construction
of a weight function as a
 certain projection of products of Drinfeld currents.
Section 4 contains the main calculations with the constructions of
\cite{EKhP} and \cite{KPT}.
For this  we extend the projection operators of \cite{ER}
to a natural completion of the algebra, where we use so called
composed currents
of \cite{DKh} and strings. Their projections are then expressed in
 Gauss coordinates of the fundamental $L$-operator.

 In Section 5 we
 translate the results
to the $L$-operator's language. The resulting formula, see Theorem 2,
 is more general then
in \cite{TV3}: it does not refer to any evaluation map. Its structure
is rather curious: the order of $L$-operators's entries does not
correspond to any normal ordering of the root system. A specialization
to the evaluation map is given by eq. \rf{Wt5}. In Section 6 we
 compare our results with  \cite{TV3}. The $R$-matrix in \cite{TV3} differs
 from ours by a finite-dimensional twist. We modify accordingly our
 construction and observe the resulting literal coincidence with \cite{TV3}. This
 allows us to verify the conjecture of \cite{KPT} for the $R$-matrix of
 \cite{TV3}: the two constructions of weight functions give the same result
 in any irreducible $\Uqgln$-module with a weight singular vector.
 This result can be formulated on a formal level, see Theorem 3.
 The paper contains three appendices with important technical results,
 including various properties of strings and of their projections.

\section{Quantum affine algebra $\Uqgln$}

\subsection{$\LL$-operator description}

\label{section2.1}

Let $\E_{ij}\in{\textrm{End}}(\CC^N)$ be a matrix with the only nonzero entry
equal to $1$
at the intersection of the $i$-th row and $j$-th column.
Let $\RR(u,v)\in{\textrm{End}}(\CC^N\ot\CC^N)\ot \CC[[{v}/{u}]]$,
\begin{equation}\label{UqglN-R}
\begin{split}
\RR(u,v)\ =\ &\ \sum_{1\leq i\leq N}\E_{ii}\ot \E_{ii}\ +\ \frac{u-v}{qu-q^{-1}v}
\sum_{1\leq i<j\leq N}(\E_{ii}\ot \E_{jj}+\E_{jj}\ot \E_{ii})
\\
+\ &\frac{q-q^{-1}}{qu-q^{-1}v}\sum_{1\leq i<j\leq N}
(v \E_{ij}\ot \E_{ji}+ u \E_{ji}\ot \E_{ij})\,,
\end{split}
\end{equation}
be a trigonometric $R$-matrix associated with the vector
representation of ${\mathfrak{gl}}_N$. It satisfies the Yang-Baxter
equation
\begin{equation}\label{YBeq}
\RR_{12}(u_1,u_2)\RR_{13}(u_1,u_3)\RR_{23}(u_2,u_3)=
\RR_{23}(u_2,u_3)\RR_{13}(u_1,u_3)\RR_{12}(u_1,u_2)\,,
\end{equation}
and the inversion relation
\begin{equation}\label{unitar}
\RR_{12}(u_1,u_2)\RR_{21}(u_2,u_1)=1.
\end{equation}

The algebra $\Uqgln$ (with the zero central charge and the gradation operator
dropped out) is a unital associative algebra generated by the modes
$\LL^{\pm}_{ij}[\pm k]$, $k\geq 0$, $1\leq i,j\leq N$, of the $\LL$-operators
$\LL^{\pm}(z)=\sum_{k=0}^\infty\sum_{i,j=1}^N \E_{ij}\otimes \LL^{\pm}_{ij}[\pm k]z^{\mp
k}$, subject to relations
\begin{equation}\label{L-op-com}
\begin{split}
&\RR(u,v)\cdot (\LL^{\pm}(u)\ot \mathbf{1})\cdot (\mathbf{1}\ot \LL^{\pm}(v))=
(\mathbf{1}\ot \LL^{\pm}(v))\cdot (\LL^{\pm}(u)\ot \mathbf{1})\cdot \RR(u,v)\\
&\RR(u,v)\cdot (\LL^{+}(u)\ot \mathbf{1})\cdot (\mathbf{1}\ot \LL^{-}(v))=
(\mathbf{1}\ot \LL^{-}(v))\cdot (\LL^{+}(u)\ot \mathbf{1})\cdot \RR(u,v),\\
&\LL^{+}_{ij}[0]=\LL^{-}_{ji}[0]=0,\qquad \LL^{+}_{kk}[0]\LL^{-}_{kk}[0]=1,
\qquad 1\leq i<j \leq N,
\quad 1\leq k\leq N\ .
\end{split}
\end{equation}

The coalgebraic structure of the algebra $\Uqgln$ is defined by the rule
\begin{equation}\label{coprL}
\Delta \sk{\LL^{\pm}_{ij}(u)}=\sum_{k=1}^N\ \LL^{\pm}_{kj}(u)\otimes
\LL^{\pm}_{ik}(u)\,.
\end{equation}

\subsection{The current realization of $\Uqgln$}
\label{section2.2}
The algebra $\Uqgln$ in the current realization (with the zero
central charge and the gradation operator dropped out) is generated
by the modes of the Cartan currents
\begin{equation}\label{current0}
k_i^\pm(z)=\sum_{m\geq0}k_i^\pm[\pm m]z^{\mp m},\qquad
k^+_i[0]k^-_i[0]=1\,, \end{equation}
 $i=1,\ldots,N$, and by the
modes of the generating functions (called 'Drinfeld
currents') 
\begin{equation}\label{currents}
E_i(z)=\sum_{n\in\ZZ}E_i[n]z^{-n}\, ,\qquad
F_i(z)=\sum_{n\in\ZZ}F_i[n]z^{-n}\, ,
\end{equation}
$i=1,\ldots,N-1$, subject to relations
$$
(q^{-1}z-q^{}w)E_{i}(z)E_{i}(w)=
E_{i}(w)E_{i}(z)(q^{}z-q^{-1}w)\, ,
$$
$$
(z-w)E_{i}(z)E_{i+1}(w)=
E_{i+1}(w)E_{i}(z)(q^{-1}z-qw)\, ,
$$
$$
(q^{}z-q^{-1}w)F_{i}(z)F_{i}(w)=
F_{i}(w)F_{i}(z)(q^{-1}z-q^{}w)\, ,
$$
$$
(q^{-1}z-qw)F_{i}(z)F_{i+1}(w)=
F_{i+1}(w)F_{i}(z)(z-w)\, ,
$$
$$
k_i^\pm(z)F_i(w)\left(k_i^\pm(z)\right)^{-1}=
\frac{q^{-1}z-qw}{z-w}F_i(w)\, ,
$$
\begin{equation}\label{gln-com}
k_{i+1}^\pm(z)F_i(w)\left(k_{i+1}^\pm(z)\right)^{-1}=
\frac{q^{}z-q^{-1}w}{z-w}F_i(w)\, ,
\end{equation}
$$
k_i^\pm(z)F_j(w)\left(k_i^\pm(z)\right)^{-1}=F_j(w)
\qquad {\rm if}\quad i\not=j,j+1\, ,
$$
$$
k_i^\pm(z)E_i(w)\left(k_i^\pm(z)\right)^{-1}=
\frac{z-w}{q^{-1}z-q^{}w}E_i(w)\, ,
$$
$$
k_{i+1}^\pm(z)E_i(w)\left(k_{i+1}^\pm(z)\right)^{-1}=
\frac{z-w}{q^{}z-q^{-1}w}E_i(w)\, ,
$$
$$
k_i^\pm(z)E_j(w)\left(k_i^\pm(z)\right)^{-1}=E_j(w)
\qquad {\rm if}\quad i\not=j,j+1\, ,
$$
$$
[E_{i}(z),F_{j}(w)]=
\delta_{{i},{j}}\ \delta(z/w)\ (q-q^{-1})\left(
k^+_{i}(z)/k^+_{i+1}(z)-k^-_{i}(w)/k^-_{i+1}(w)\right)\, ,
$$
together with the Serre relations
\begin{equation}
\begin{split}
{\rm Sym}_{z_1,z_{2}}
(E_{i}(z_1)E_{i}(z_2)E_{i\pm 1}(w)
&-(q+q^{-1})E_{i}(z_1)E_{i\pm 1}(w)E_{i}(z_2)+\\
&+E_{i\pm 1}(w)E_{i}(z_1)E_{i}(z_2))=0\, ,\\
\label{serre}
{\rm Sym}_{z_1,z_{2}}
(F_{i}(z_1)F_{i}(z_2)F_{i\pm 1}(w)
&-(q+q^{-1})F_{i}(z_1)F_{i\pm 1}(w)F_{i}(z_2)+\\
&+F_{i\pm 1}(w)F_{i}(z_1)F_{i}(z_2))=0\, .
\end{split}
\end{equation}
The isomorphism of the two realization is established with a help of
the Gauss decomposition of $\LL$-operators. We define Gauss coordinates
$\FF^{\pm}_{j,i}(z)$, $\EE^{\pm}_{i,j}(z)$ $1\leq i<j\leq N$ and
$k^\pm_{i}(z)$, $i=1,...,N$ by the relations
\begin{equation}\label{L-op}
\LL^\pm(z)=\sk{\sum_{i=1}^N
\E_{ii}+\sum^N_{i<j}\FF^{\pm}_{j,i}(z)\E_{ij}}
\cdot\sum^N_{i=1}k^\pm_{i}(z)\E_{ii}\cdot \sk{\sum_{i=1}^N
\E_{ii}+\sum^N_{i<j}\EE^{\pm}_{i,j}(z)\E_{ji}}\,.
\end{equation}
We identify $k_i^\pm(z)$ with the Cartan currents \rf{current0} and set
\cite{DF}
\begin{equation}\label{DF-iso}
E_i(z)=\EE^{+}_{i,i+1}(z)-\EE^{-}_{i,i+1}(z)\,,\quad
F_i(z)=\FF^{+}_{i+1,i}(z)-\FF^{-}_{i+1,i}(z)\ .
\end{equation}

{}For any series $\acc(t)=\sum_{m\in\ZZ}\acc[m]t^{-m}$ we denote
$\acc(t)^{(+)}=\sum_{m>0} \acc[m]\, t^{-m}\,,
$ and $\acc(t)^{(-)}=-\sum_{ m\leq 0} \acc[m]\, t^{-m}\,.\
$
The initial conditions \rf{L-op-com} imply the relations
\begin{equation}\label{DFinverse}
\FF^{\pm}_{i+1,i}(z)=F_i(z)^{(\pm)},\qquad
\EE^{\pm}_{i,i+1}(z)=z\left(z^{-1}E_i(z)\right)^{(\pm)}.
\end{equation}

 In \cite{D88} the \emph{current} Hopf structure for the algebra
$\Uqgln$ has been defined,
\begin{equation}\label{gln-copr}
\begin{split}
\Delta^{(D)}\sk{E_i(z)}&=E_i(z)\ot 1 + k^-_{i}(z)\sk{k^-_{i+1}(z)}^{-1}\ot
E_i(z),\\
\Delta^{(D)}\sk{F_i(z)}&=1\ot F_i(z) + F_i(z)\ot
k^+_{i}(z)\sk{k^+_{i+1}(z)}^{-1},\\
\Delta^{(D)}\sk{k^\pm_i(z)}&=k^\pm_i(z)\ot k^\pm_{i}(z).
\end{split}
\end{equation}

We consider two types of Borel subalgebras of the
algebra $\Uqgln$. Borel subalgebras $U_q(\mathfrak{b}^\pm)\subset \Uqgln $
are generated by the modes of the $\LL$-operators $\LL^{(\pm)}(z)$, respectively.
Another type of Borel subalgebras is related to the current realization of
$\Uqgln$. The Borel subalgebra $U_F\subset \Uqgln$ is generated by
the modes
$F_i[n]$, $k^+_j[m]$, $i=1,\ldots,N-1$, $j=1,\ldots,N$, $n\in\ZZ$
and $m\geq0$. The Borel subalgebra
$U_E\subset \Uqgln$ is generated by the modes
$E_i[n]$, $k^-_j[-m]$, $i=1,\ldots,N-1$, $j=1,\ldots,N$, $n\in\ZZ$ and
$m\geq0$. We  also consider  a subalgebra $U'_F\subset U_F$, generated
by the elements
$F_i[n]$, $k^+_j[m]$, $i=1,\ldots,N-1$, $j=1,\ldots,N$, $n\in\ZZ$
and $m>0$, and a subalgebra $U'_E\subset U_E$ generated by
the elements
$E_i[n]$, $k^-_j[-m]$, $i=1,\ldots,N-1$, $j=1,\ldots,N$, $n\in\ZZ$
and $m>0$.
Further, we will be interested in the intersections,
\begin{equation}
\label{Intergl}
U_f^-=U'_F\cap U_q(\mathfrak{b}^-)\,,\qquad
U_F^+=U_F\cap U_q(\mathfrak{b}^+)\,
\end{equation}
and will describe properties of projections to these intersections.

It was proved in \cite{KPT} that the subalgebras $U_f^-$
and $U_F^+$ are coideals with respect to Drinfeld coproduct
\r{gln-copr}
\begin{equation*}
\Delta^{(D)}(U_F^+)\subset \Uqgln\ot U_F^+\,,\qquad
\Delta^{(D)}(U_f^-)\subset U_f^-\ot \Uqgln\,,
\end{equation*}
and the multiplication $m$ in $\Uqgln$ induces an isomorphism of vector spaces
$$m: U_f^-\ot U_F^+\to U_F\,.$$
According to the general theory presented in \cite{EKhP} we  define
 projection operators $\Pfp:U_F\subset \Uqgln \to U_F^+$ and
$\Pfm:U_F\subset \Uqgln \to U_f^-$ by the prescriptions
\begin{equation}\label{pgln}
\begin{split}
\Pfp(\Ff_-\ \Ff_+)&=\coun(\Ff_-)\ \Ff_+, \qquad
\Pfm(\Ff_-\ \Ff_+)=\Ff_-\ \coun(\Ff_+),
\\& \text{for any}\qquad \Ff_-\in U_f^-,
\quad \Ff_+\in U_F^+ .
\end{split}
\end{equation}

Denote by  $\overline U_F$ an extension of the algebra $U_F$ formed
by linear combinations of series, given as infinite sums of monomials
$a_{i_1}[n_1]\cdots a_{i_k}[n_k]$ with $n_1\leq\cdots\leq n_k$, and $n_1+...+n_k$
fixed,
where  $a_{i_l}[n_l]$ is either $F_{i_l}[n_l]$ or $k^+_{i_l}[n_l]$.
It was proved in \cite{EKhP} that
\begin{itemize}
\item[(1)] the action of the projections \r{pgln} can be extended to the
 agebra
$\overline U_F$;
\item[(2)] for any $\Ff\in \overline U_F$ with $\Delta^{(D)}(\Ff)=\sum_i \Ff'_i\otimes \Ff''_i$ we have
\begin{equation}\label{pr-prop}
\Ff=\sum_i\Pfm(\Ff''_i)\cdot \Pfp(\Ff'_i)\,.
\end{equation}
\end{itemize}

\section{Weight functions}
\subsection{Definitions}
\label{section3}
 This section is based on the paper
\cite{KPT}\footnote{The definition of the weight function used in
this paper differs from that of \cite{KPT} by the reverse
 order of the variables. The definition of the modified weight
function is the same.}.

We call a vector $v$ \emph{a weight singular vector} if it is
annihilated by any non-negative mode  $E_i[n]$, $i=1,\ldots,N-1$,
$n\geq 0$ and is an eigenvector for  $k^+_i(z)$, $i=1,\ldots,N$
\begin{equation}\label{hwv}
\EE^{+}_{i,i+1}(z)\cdot v=0\ ,\qquad k^+_i(z)\cdot v=\lambda_i(z)\,v\,,
\end{equation}
where $\Lambda_i(z)$ is a meromorphic function, decomposed as a
power series in $z^{-1}$. The $L$-operator \r{L-op}, acting on a
weight singular vector $v$, becomes upper-triangular
\begin{equation}\label{L-op-tr}
\LL^{+}_{ij}(z)\ v=0\,,\quad i>j\,,\quad \LL^{+}_{ii}(z)\ v=\lambda_i(z)\ v\,,\quad
i=1,\ldots,N\,.
\end{equation}

We define a weight function by its comultiplication properties.

Let $\Pi$ be the set $\{1,\ldots,N-1\}$ of indices of the simple
positive roots of $\mathfrak{gl}_N$. A finite collection
$I=\{i_1,\dots,i_n\}$ with a linear ordering $i_i\ord\cdots\ord i_n$
and a `coloring' map $\iota:I\to\Pi$ is called an {\it ordered
$\Pi$-multiset}. Sometimes, we denote the map $\iota$ by $\iota_I$.
A morphism between two ordered $\Pi$-multisets $I$ and $J$ is a map
$m:I\to J$ that respects the orderings in $I$ and $J$ and
intertwines the maps $\iota_I$ and $\iota_J$: $\iota_J m=m\iota_I$.
In particular, any subset $I'\subset I$ of a $\Pi$-ordered multiset
has a unique structure of $\Pi$-ordered multiset, such that the
inclusion map is a morphism of $\Pi$-ordered multisets.

To each $\Pi$-ordered multiset $I=\{i_1,\dots,i_n\}$ we attach an
ordered set of variables $\{t_i|i\in I\}=\{t_{i_1},\dots,t_{i_n}\}$.
Each element $i_k\in I$ and each variable $t_{i_k}$ has its own
`type': $\iota(i_k)\in\Pi$.

Let $i$ and $j$ be elements of some ordered $\Pi$-multiset. Define a
rational function
\begin{equation}\label{gam}
\gamma(t_i,t_j)=\left\{
\begin{array}{ll}
\ds\frac{t_i-t_j}{qt_i-q^{-1}t_j}\ ,\quad&\mbox{if}\quad \i(i)=\i(j)+1\ ,\\[5mm]
\ds\frac{q^{-1}t_i-qt_j}{t_i-t_j}\ ,\quad&\mbox{if}\quad \i(j)=\i(i)+1\ ,\\[5mm]
\ds\frac{q^{}t_i-q^{-1}t_j}{q^{-1}t_i-q^{}t_j}\ ,\quad&\mbox{if}\quad \i(i)=\i(j)\ ,\\
[5mm] 1\ ,\quad&\mbox{otherwise}\ .
\end{array}
\right.
\end{equation}

Assume that for any representation $V$ of $\Uqgln$ with a weight
singular vector $v$, and any ordered $\Pi$-multiset $I=\{i_1,\dots i_n\}$,
there is a $V$-valued rational function $w_{V,I}(t_{i_1},\ldots,t_{i_n})\in V$
depending on the variables $\{t_i|i\in I\}$. We call such a collection of
rational functions a {\it weight function} $w$, if:
\begin{itemize}
\item[(a)] The rational function, corresponding to the empty set,
is equal to $v$:
$w_{V,\emptyset}\equiv v$\, .
\item[(b)] The function $w_{V,I}(t_{i_1},\ldots,t_{i_n})$ depends only on
an isomorphism class of an ordered $\Pi$-multiset, that is, for any isomorphism
$f:I\to J$ of ordered $\Pi$-multisets
\begin{equation}
\label{isomorph}
w_{V,I}(t_{f(i)}|_{i\in I})=w_{V,J}(t_j|_{j\in J})\,.
\end{equation}
\item[(c)]
The functions $w_{V,I}$ satisfy the following comultiplication
property. Let $V=V_1\otimes V_2$ be a tensor product of two
representations with singular vectors $v_1$, $v_2$ and weight series
$\{\lambda_b^{(1)}(u)\}$ and $\{\lambda_b^{(2)}(u)\}$,
$b=1,\ldots,N$. Then for any ordered $\Pi$-multiset $I$ we have
\end{itemize}
\begin{equation}\label{weight2}
w_{V,I}(t_{i} |_{i\in I})=
\sum\limits_{I=I_1\coprod I_2}
w_{V_1,I_1}(t_{i}|_{i\in I_1})\otimes w_{V_2,I_2}(t_{i}|_{i\in I_2})
\ \Phi_{I_1,I_2}(t_{i}|_{i\in I})
\prod\limits_{{j\in I_1}}
\frac{\lambda^{(2)}_{\iota(j)}(t_{j})}{\lambda^{(2)}_{\iota(j)+1}(t_{j})}
\,,
\end{equation}
{where}
\begin{equation}
\Phi_{I_1,I_2}(t_{i}|_{i\in I})=
\prod\limits_{
{i\in I_2,\ j\in I_1,\ i\ord j}} \gamma(t_i,t_j).
\end{equation}
The summation in \r{weight2} runs over all possible decompositions of
the set $I$ into a disjoint union of two non-intersecting
subsets $I_1$ and $I_2$. The structure of ordered $\Pi$-multiset on each subset
is induced from that of $I$.

Given elements $i,j$ of some ordered multiset define two functions
$\tilde{\gamma}(t_i,t_j)$ and $\beta(t_i,t_j)$ by the formulae
\begin{equation*}
\tilde\gamma(t_i,t_j)=\left\{
\begin{array}{ll}
\ds\frac{t_i-t_j}{qt_i-q^{-1}t_j}\ ,\quad&\mbox{if}\quad \i(i)=\i(j)+1\ ,\\[5mm]
\ds\frac{q^{-1}t_i-qt_j}{t_i-t_j}\ ,\quad&\mbox{if}\quad \i(j)=\i(i)+1\ ,\\[5mm]
1\ ,\quad&\mbox{otherwise}
\end{array}
\right.
\end{equation*}
and
\begin{equation}\label{beta}
\beta(t_i,t_j)=\left\{
\begin{array}{ll}
\ds\frac{q^{-1}t_i-qt_j}{t_i-t_j}\ ,\quad&\mbox{if}\quad \i(i)=\i(j)\ ,\\[5mm]
1\ ,\quad&\mbox{otherwise\,.}
\end{array}
\right.
\end{equation}

A collection of rational $V$-valued functions $\w_{V,I}(t_{i}|_{i\in
I})$, depending on a representation $V$ of $\Uqgln$ with a
weight singular vector $v$, and an ordered $\Pi$-multiset $I$, is
called a {\it modified weight function} $\w$, if it satisfies
conditions (a), (b) above  and the condition (c$'$):

\begin{itemize}
\item[(c$'$)]
Let $V=V_1\otimes V_2$ be a tensor product of two representations
with singular vectors $v_1$, $v_2$ and weight series
$\{\lambda_b^{(1)}(u)\}$ and $\{\lambda_b^{(2)}(u)\}$,
$b=1,\ldots,N$. Then for any multiset $I$ we have
\end{itemize}
\begin{equation}\label{weight22}\begin{split}
\w_{V,I}(t_{i} |_{i\in I})
=&\sum\limits_{I=I_1\coprod I_2}
\w_{V_1,I_1}(t_{i}|_{i\in I_1})
\otimes
\w_{V_2,I_2}(t_{i}|_{i\in I_2})
\cdot\tilde{\Phi}_{I_1,I_2}(t_{i}|_{i\in I})\times\\
\times&\prod\limits_{{j\in I_1}}
{\lambda^{(2)}_{\iota(j)}(t_{j})}
\prod\limits_{{j\in I_2}}
{\lambda^{(1)}_{\iota(j)+1}(t_{j})},
\end{split}\end{equation}
where
\begin{eqnarray*}
\tilde{\Phi}_{I_1,I_2}(t_{i}|_{i\in I})= \prod_{i\in I_1,\ j\in I_2}
\beta(t_i,t_j) \prod_{{i\in I_2,\ j\in I_1,\ i\ord
j}}\tilde\gamma(t_i,t_j)\,.
\end{eqnarray*}

There is a bijection between weight functions
and modified weight functions given by the following relations.
Let $w$ be a weight function. Then the collection
$\w_{V,I}(t_i|_{i\in I})$, where
\begin{equation}\label{swf}
\w_{V,I}(t_i|_{i\in I})=w_{V, I}(t_i|_{i\in I})
\prod\limits_{i\ord j } \beta(t_i,t_j)\,
\prod_{i\in I}
\lambda_{\i(i)+1}(t_i)
\end{equation}
is a modified weight function.
%

\subsection{Weight function and projections}

Let $I=\{i_1,\dots,i_n\}$ and $J=\{j_1,\dots,j_n\}$ be two ordered
$\Pi$-multisets. Let $\sigma:I\to J$ be an invertible map, which
intertwines the coloring maps $\iota_I$ and $\iota_J$: $\iota_J
\sigma=\sigma\iota_I$, but does not necessarily respect the
orderings in $I$ and $J$ (that is, $\sigma$ is a `permutation' on
classes of isomorphisms of ordered $\Pi$-multisets).

Let $w(t_j|_{j\in J})$ be a function of the variables $t_j|_{j\in J}$.
Define a {\it pullback} $^{\sigma,\gamma}\!w(t_i|_{i\in I})$ by the rule
\begin{equation}\label{pullback}
^{\sigma,\gamma}\!w(t_i|_{i\in I})=w(t_{\sigma(i)}|_{i\in I})
\prod_{
i,j\in I,\,
i\ord j,\ \sigma(j)\ord \sigma(i)
}\gamma(t_i,t_j)\,.
\end{equation}

One may check \cite{KPT} that the pullback operation \rf{pullback}
is compatible with the comultiplication rule \rf{weight2}. We call a
weight function $w_{V,I}(t_i|_{i\in I})$ {\it $q$-symmetric}, if for
any ordered $\Pi$-multisets $I$ and $J$ and an invertible map
$\sigma:I\to J$, intertwining the colouring maps, we have
\begin{equation}\label{q-s}
^{\sigma,\gamma}w_{V,J}(t_{i} |_{i\in I})= w_{V,I}(t_{i} |_{i\in I})\,.
\end{equation}

For any multiset $I$ we define a $U^+_F$-valued series $\mathcal{W}_I(t_i|_{i\in I})$
as the projection
\begin{equation}\label{W-un-any}
\mathcal{W}_I(t_i|_{i\in I})=\Pfp\sk{F_{\i(i_n)}(t_{i_n})
F_{\i(i_{n-1})}(t_{i_{n-1}})\cdots F_{\i(i_2)}(t_{i_2}) F_{\i(i_1)}(t_{i_1})}
\end{equation}
We call the  series  $\mathcal{W}_{I}(t_i|_{i\in I})$ {\em  universal
weight function}.

\bigskip

\noindent {\bf Theorem}\ \ \cite{KPT}. {\em A collection of\/ $V$-valued functions
 $w_{V,I}(t_i|_{i\in I})=\mathcal{W}_I(t_i|_{i\in I})\ v$, where $V$ is a
$\Uqgln$-module with a singular weight vector $v$,
form a $q$-symmetric weight function.}

\bigskip

 Let $\bar n=\{n_{1},n_{2},\ldots,n_{N-2},n_{N-1}\}$
be a set of non-negative integers. Let $I_{{\bar n}}$, be an ordered $\Pi$-multiset,
 such that its first $n_{1}$ elements have the type $1$, the next
$n_2$ elements have the type 2 and the last
$n_{N-1}$ elements have the type $N-1$. Denote by  $\bar t_{\segg{\bar n}}$
 the set of related variables:
\begin{equation}\label{set111}
\bar{t}_{\segg{\bar{n}}} = \left\{ t^{1}_{1},\ldots,t^{1}_{n_{1}};
t^{2}_{1},\ldots, t^{2}_{n_{2}};\ \ldots\ldots\ ;
t^{N-2}_{1},\ldots, t^{N-2}_{n_{N-2}};
t^1_{N-1},\ldots,t^{N-1}_{n_{N-1}}\right\}\,.
\end{equation}
The variable $t_{k}^a$ is of type $a$. If $n_a=0$ for some $a$, then
the variables of the type $a$ are absent in the set \r{set111}.
Denote by $\mathcal{W}^{N-1}(\bar{t}_{\segg{\bar{n}}})$
the universal weight function associated with the set of  variables
\r{set111}.
\begin{equation}\label{uwf}
\mathcal{W}^{N-1}(\bar{t}_{\segg{\bar{n}}})=P^+\left(
F_{N-1}(t_{n_{N-1}}^{N-1})\cdots F_{N-1}(t_{1}^{N-1})
  \quad \cdots\quad
F_{1}(t_{n_1}^1)\cdots F_{1}(t_{1}^1)
\right).
\end{equation}
 For any weight singular vector $v$ let ${w}_V^{N-1}(\bar
t_{\segg{\bar n}})=\mathcal{W}^{N-1}(\bar t_{\segg{\bar n}})\ v$ be
the related weight function and
\begin{equation}\label{b-v-sym}
{\bf w}_V^{N-1}(\bar t_{\segg{\bar n}})=\qsym(\bar t_{\segg{\bar
n}})\prod_{a=2}^{N}\prod_{\ell=1}^{n_{a-1}}\lambda_a(t^{a-1}_\ell) \
{w}^{N-1}(\bar t_{\segg{\bar n}})\end{equation} the corresponding modified
weight function. Here
\begin{equation*}
\qsym(\bar t_{\segg{\bar n}})=\prod_{a=1}^{N-1}\prod_{1\leq \ell<\ell'\leq n_a}
\frac{q-q^{-1}t^a_\ell/t^a_{\ell'}}{1-t^a_\ell/t^a_{\ell'}}\ .
\end{equation*}
We call the modified weight function \r{b-v-sym} {\em  off-shell Bethe vector}.

One can see that any $\Pi$-ordered multiset is isomorphic to a permutation of some
$I_{\bar{n}}$.  The $q$-symmetric property then implies that the
series \rf{uwf} completely describe the universal weight function
\rf{W-un-any}, and of-shell Bethe vectors \r{b-v-sym} completely
describe the corresponding modified weight function.
\section{A computation of the universal weight function}\label{sec4}

In this section we express the universal weight function
$\mathcal{W}^{N-1}(\bar{t}_{\bar{n}})$ in generators of $\Uqgln$,
using the definition of the projection operator $P^+$. Our  strategy
is as follows. Under projection operator in \r{uwf} we separate all
factors $F_a(t^a_i)$ with $a<N-1$ and apply to this product the
ordering procedure of Proposition \ref{decomp}, based on the
property \rf{dec-ff1}. We get under total projection a
symmetrization of a sum of terms $x_iP^-(y_i)P^+(z_i)$ with rational
functional coefficients; here each $x_i$ is a monomial on the modes of
$F_{N-1}(t)$, and $y_i$, $z_i$ are monomials on the modes of $F_a(t)$ with $a<N-1$.
Then we reorder $x_i$ and $P_-(y_i)$. At this stage  composed
currents, collected in so called strings, appear. The calculation of
the projection of strings is a separate problem, which is solved by
analytical tools.

\subsection{Basic notations}

Let $\bar\lll$ and $\bar\rr$ be two collections of nonnegative integers
satisfying a set of inequalities
\begin{equation}\label{set23}
\lll_a\leq\rr_a\,,\quad a=1,\ldots,N-1\,.
\end{equation}
Denote by $\seg{\bar\rr}{\bar\lll}$ a collection of segments which
are sets of positive increasing integers
$\{\lll_a+1,\ldots,\rr_a-1,\rr_a\}$ including $\rr_a$ and excluding
$\lll_a$. The length of each segment is equal to $\rr_a-\lll_a$.

For a given set $\seg{\bar\rr}{\bar\lll}$ of segments  we denote by
$\bar t_{\seg{\bar\rr}{\bar\lll}}$ the set of variables
\begin{equation}\label{set2}
\bar{t}_{\seg{\bar{\rr}}{\bar\lll}}\ =
\{t^{1}_{\lll_{1}+1},\ldots,t^{1}_{\rr_{1}} ;\ \ldots\ ;
t^{N-2}_{\lll_{N-2}+1},\ldots,t^{N-2}_{\rr_{N-2}};\,
t^{N-1}_{\lll_{N-1}+1},\ldots,t^{N-1}_{\rr_{N-1}}\}\,.\end{equation}
 The
number of the variables of the type $a$ is equal to $\rr_a-\lll_a$.
In this notation, the set of the variables \r{set111} is $\bar
t_{\segg{\bar n}}\equiv \bar t_{\seg{\bar n}{\bar 0}}$. One can
consider \rf{set2} as a list of variables, corresponding to the
ordered multiset, naturally related to $\seg{\bar\rr}{\bar\lll}$.

{}For any $a=1,\ldots,N-1\,$ we consider the segment
$\seg{\rr_a}{\lll_a}=\{\lll_a+1,\ldots,\rr_a-1,\rr_a\}$
as an ordered multiset $\{\lll_a+1\prec\cdots\prec\rr_a-1\prec\rr_a\}$,
in which all the elements are of the type $a$. The related set of
variables is denoted as
\begin{equation}\label{set3}
\bar{t}^a_{\seg{{\rr_a}}{\lll_a}}\
=\{
t^{a}_{\lll_{a}+1},\ldots,t^{a}_{\rr_{a}} \}\,.\end{equation}
All the variables in \rf{set3} have  the type $a$. For the segment
$\seg{\rr_a}{\lll_a}=\seg{n_a}{0}$ we use the shorten notation
$\bar{t}^a_{\seg{{n_a}}{0}}\equiv\bar{t}^a_{\segg{n_a}}$.

Our basic calculations are performed on a level of formal series attached
 to certain ordered multisets. For the save of space we often write some series
as rational homogeneous functions with the following prescription.
Let  $\{t_i|i\in I\}=\{t_{i_1},\dots,t_{i_n}\}$ be the ordered set of
variables attached to an ordered set $I=\{i_1\prec
i_2\prec\cdots\prec i_n\}$ and $g(t_i|i\in I)$ be a rational
function. Then we associate to $g(t_i|i\in I)$ a Loran series which
is the expansion of $g(t_i|i\in I)$ in a region
 $|t_{i_1}| \ll |t_{i_2}| \ll \cdots \ll |t_{i_n}|$. If, for instance,
$1\prec 2$, then we associate to a rational function $\dfrac{1}{t_1-t_2}$ a series
$-\sum_{k\geq 0}t_1^k t_2^{-k-1}$.  With this convention to a rational function
of the variables $\bar{t}_{\segg{\bar{n}}}$ we associate a Taylor series on
$t^b_k/t^c_l$ with $b<c$ and on $t^a_i/t^a_j$ with $i<j$.

For a collection of variables $\bar t_{\seg{\bar\rr}{\bar\lll}}$ we
consider an ordered product
\begin{equation}\label{FFFl}
\F(\bar t_{\seg{\bar\rr}{\bar\lll}})=\!\!\prod_{N-1\geq a\geq
1}^{\longleftarrow} \sk{\prod_{\rr_a\geq\ell>
\lll_a}^{\longleftarrow} F_{a}(t^a_{\ell})}=
F_{N-1}(t^{N-1}_{\rr_{N-1}})\cdots
F_{1}(t^{1}_{\rr_{1}})\cdots F_{1}(t^1_{\lll_1+1}),
\end{equation}
where the series $\,F_a(t)\equiv F_{a+1,a}(t)\,$ is defined by
\r{currents}. As a particular case, we have
$\,\F(\bar t^a_{\seg{\bar\rr_a}{\bar\lll_a}})=
F_{a}(t^a_{r_a})\cdots F_{a}(t^a_{l_a+2})
F_{a}(t^a_{l_a+1})$.

 Symbols $\mathop{\prod}\limits^{\longleftarrow}_a A_a$ and
$\mathop{\prod}\limits^{\longrightarrow}_a A_a$
 mean  ordered products of
noncommutative entries $A_a$, such that $A_a$ is on the right (resp., on the left)
 from $A_b$ for $b>a$:
\begin{equation*}
\mathop{\prod}\limits^{\longleftarrow}_{j\geq a\geq i} A_a = A_j\,A_{j-1}\,
\cdots\, A_{i+1}\,A_i\,,\quad
\mathop{\prod}\limits^{\longrightarrow}_{i\leq a\leq j} A_a = A_i\,A_{i+1}\,
\cdots\, A_{j-1}\,A_j
\end{equation*}

For  collections of positive integers $\bar\rr$ and $\bar\lll$ which
satisfy inequalities \r{set23} we denote by
$S_{{\bar\lll,\bar\rr }} =S_{\rr_{N-1}-\lll_{N-1}}\times \cdots
\times S_{\rr_1-\lll_1}$ a direct product of the symmetric groups.
The group $S_{{\bar\lll,\bar\rr}}$ naturally acts on functions
of the variables $\bar t_{\seg{\bar\rr}{\bar\lll}}$ by permutations
of variables of the same type. If $\si=\si^{N-1}\times\cdots\times
\si^1\in S_{{\bar\lll,\bar\rr}}$, set
\begin{equation}\label{sigmat}
^\si \bar t_{\seg{\bar\rr}{\bar\lll}} =\{t^{N-1}_{\si^{N-1}(\lll_{N-1}+1)},\ldots,
t^{N-1}_{\si^{N-1}(\rr_{N-1})};
\ldots;t^1_{\si^1(\lll_1+1)},\ldots,t^1_{\si^1(\rr_1)}\}.
\end{equation}

For a formal series or a function $G(\bar t_{\bar\rr,\bar\lll})$ the
$q$-symmetrization means
\begin{equation}\label{qs1r}
\qSym_{\ \bar t_{\seg{\bar\rr}{\bar\lll}}} \ G(\bar
t_{\seg{\bar\rr}{\bar\lll}})= \sum_{\si\in
S_{\bar\lll,\bar\rr}}\prod_{1\leq a\leq N-1}
\prod_{\substack{\ell>\ell'\\ \si^a(\ell)<\si^a(\ell')}}
\frac{q^{}-q^{-1}t^a_{\si^a(\ell)}/t^a_{\si^a(\ell')}}
{q^{-1}-q^{}t^a_{\si^a(\ell)}/t^a_{\si^a(\ell')}}\ G(^\si \bar
t_{\seg{\bar\rr}{\bar\lll}})\,.
\end{equation}
The operation \r{qs1r} is well-defined throughout the paper, see \cite{KP} for details.
According to \r{q-s} we call a formal series $G(\bar t_{\seg{\bar\rr}{\bar\lll}})$
$q$-symmetric if
\begin{equation}\label{exa3}
\qSym%
_{\ \bar t_{\seg{\bar\rr}{\bar\lll}}}
\sk{G(\bar t_{\seg{\bar\rr}{\bar\lll}})}=\prod_{a=1}^{N-1} (\rr_a-\lll_a)!\
G(\bar t_{\seg{\bar\rr}{\bar\lll}})\,.
\end{equation}
The $q$-symmetrization of any series is  $q$-symmetric.

For a set of segments $\seg{\bar\rr}{\bar\lll}$ we introduce a third collection of
nonnegative integers $\bar\ss$ such that
 $0\leq \ss_a\leq \rr_a-\lll_a$, $a=1,\ldots,N-1$.
Each integer $\ss_a$ divides the segment $\seg{\bar\rr}{\bar\lll}$ into
two nonintersecting
segments $\seg{\bar\rr-\bar\ss}{\bar\lll}$ and
$\seg{\bar\rr}{\bar\rr-\bar\ss}$. Recall that we include
into segment its left edge and exclude the right one. It means that
$\seg{\rr_a-\ss_a}{\lll_a}=\{\lll_a+1,\ldots,\rr_a-\ss_a-1,\rr_a-\ss_a\}$ and
$\seg{\rr_a}{\rr_a-\ss_a}=\{\rr_a-\ss_a+1,\ldots,\rr_a-1,\rr_a\}$.

For a set of the variables $\bar t_{\seg{\bar\rr}{\bar\lll}}$ and a collection of integers
$\bar\ss$ which divide the set of segments $\seg{\bar\rr}{\bar\lll}$ we define a series
\begin{equation}\label{Zser}
Z_{\bar\ss}({\bar t}_{\seg{\bar\rr}{\bar\lll}})=\prod_{a=1}^{N-2}\ \
\prod_{\substack{ \rr_a-\ss_a < \ell\leq \rr_a \\
  \lll_{a+1}< \ell' \leq \rr_{a+1}-\ss_{a+1}}} \frac{q-q^{-1}\
t^{a}_{\ell}\,/\,t^{a+1}_{\ell'}}{1-t^{a}_{\ell}\,/\,t^{a+1}_{\ell'}}\,.
\end{equation}
Note that this series does not depend on the variables
$t^{N-1}_{\rr_{N-1}-\ss_{N-1}+1},\ldots,t^{N-1}_{\rr_{N-1}}$ and
$t^1_{\lll_1+1},\ldots,t^1_{\rr_1-\ss_1}$. If $\ss_{N-1}=\rr_{N-1}$
then this series does not depend
on all variables of the type $N-1$ from the set $ t_{\seg{\bar\rr}{\bar\lll}}$.
Also if $\rr_a=\lll_a$ for all $a$ except one value $a=j$  then
 the collection $\seg{\bar\rr}{\bar\lll}$ of segments contains only one segment
of the type $j$ and  we set the  series \r{Zser} equal to 1.

We call any expression $\sum_i f^{(i)}_-\cdot f^{(i)}_+$, where
$f^{(i)}_-\in U^-_f$ and $f^{(i)}_+\in U^+_f$ {\it (normal) ordered}.
Using the property \r{pr-prop} of the projections
we can present any product \r{FFFl} in a normal ordered form.
\begin{proposition}\label{decomp}We have an equality
\begin{equation}\label{dec-ff1}
\begin{array}{c}
\ds \F(\bar t_{\seg{\bar\rr}{\bar\lll}})=
\sum_{0\leq \ss_{N-1}\leq \rr_{N-1}-\lll_{N-1}}\cdots \sum_{0\leq \ss_1\leq \rr_1-\lll_1}
\ \  \prod_{1\leq a\leq N-1}  \frac{1}{(\ss_a)!(\rr_a-\lll_a-\ss_a)!}\times \\ [10mm]
\ds \times\   \qSym_{\ \bar t_{\seg{\bar\rr}{\bar\lll}}}
\left(Z_{\bar\ss}({\bar t}_{\seg{\bar\rr}{\bar\lll}})
\ds\  \Pfm\sk{\F(\bar t_{\seg{\bar\rr}{\bar\rr-\ss}})}\cdot
\Pfp\sk{\F(\bar t_{\seg{\bar\rr-\bar\ss}{\bar\lll}})}\right).
\end{array}
\end{equation}
\end{proposition}
\noindent {\it Proof}. We use the coproduct \r{gln-copr}. The $q$-symmetrization appears
in \r{dec-ff1}  due
to the commuting of the Cartan and simple root currents corresponding to
the same type, and the
series $Z_{\bar s}({\bar t}_{\seg{\bar\rr}{\bar\lll}})$ appears due to the relation
\begin{eqnarray*}
&k^+_{a+1}(t^{a+1}_{\rr_{a+1}-\ss_{a+1}})\cdots k^+_{a+1}(t^{a+1}_{\lll_{a+1}+1})\cdot
F_a(t^a_{\rr_a})\cdots F_a(t^a_{\rr_a-\ss_a+1})=\\
&=\ds \prod_{\substack{ \rr_a-\ss_a < \ell\leq \rr_a \\
  \lll_{a+1}< \ell' \leq \rr_{a+1}-\ss_{a+1}}}
 \frac{q-q^{-1}\
t^{a}_{\ell}\,/\,t^{a+1}_{\ell'}}{1-t^{a}_{\ell}\,/\,t^{a+1}_{\ell'}}\
F_a(t^a_{\rr_a})\cdots F_a(t^a_{\rr_a-\ss_a+1})\times \\
&\quad\times
k^+_{a+1}(t^{a+1}_{\rr_{a+1}-\ss_{a+1}})\cdots k^+_{a+1}(t^{a+1}_{\lll_{a+1}+1})
\end{eqnarray*}
which has to be used in order to apply the operator $\Pfm$ in the right hand side of
\r{dec-ff1}.

\subsection{Composed currents and strings}
\label{MoreDF}

Following \cite{DKh,KP}, we introduce  {\it composed currents\/}
$\ff_{j,i}(t)$ for $i<j$. The composed currents for
 nontwisted quantum affine algebras were defined in \cite{DKh}.
The series $\ff_{i+1,i}(t)$, $i=1,\ldots N-1$,  coincides with
$\ff_i(t)$, cf.~\r{currents}.
According to \cite{DKh}, the coefficients of the series $\ff_{j,i}(t)$ belong to
the completion $\overline U_F$ of the algebra $U_F$, see Section \ref{section2.2}.

The completion $\overline U_F$ determines analyticity properties of products
of currents (and coincide with analytical properties of their matrix coefficients for highest
weight representations). One can show that for $|i-j|>1$, the product
$F_i(t)F_j(w)$ is an expansion of a function analytic at $t\ne 0$, $w\ne 0$.
The situation is more delicate for $j=i,i\pm1$. The products $F_i(t)F_i(w)$ and
$F_i(t)F_{i+1}(w)$ are expansions of analytic functions at $|w|<|q^2 t|$, while
the product $F_i(t)F_{i-1}(w)$ is an expansion of an analytic function at
$|w|<|t|$. Moreover, the only singularity of the corresponding functions
in the whole region $t\ne 0$, $w\ne 0$, are simple poles at the respective
hyperplanes, $w=q^2t$ for $j=i,i+1$, and $w=t$ for $j=i-1$.

The definition of the composed currents may be written in analytical form
\begin{equation}\label{rec-f1}
\ff_{j,i}(t)\, =\,
-\,\mathop{\rm res}\limits_{w=t}\ff_{j,a}(t)\ff_{a,i}(w)\,\frac{dw}{w}\, =\,
\mathop{\rm res}\limits_{w=t}\ff_{j,a}(w)\ff_{a,i}(t)\,\frac{dw}{w}
\end{equation}
for any $a=i+1,\ldots, j-1$. It is equivalent to the relation
\begin{equation}\label{rec-f111}
\begin{split}
\ff_{j,i}(t)&=
\oint \ff_{j,a}(t) \ff_{a,i}(w)\ \frac{dw}{w}-
\oint \frac{q^{-1}-qt/w}{1-t/w}
\;\ff_{a,i}(w) \ff_{j,a}(t)\ \frac{dw}{w}\,,\\
\ff_{j,i}(t)&=
\oint \ff_{j,a}(w) \ff_{a,i}(t)\ \frac{dw}{w}-
\oint  \frac{q^{-1}-qw/t}{1-w/t}
\;\ff_{a,i}(t) \ff_{j,a}(w)\ \frac{dw}{w}\,.\\
\end{split}
\end{equation}
In \r{rec-f111} $\oint \frac{dw}{w} g(w)=g_0$ for any formal series $g(w)=\sum_{n\in\ZZ}g_n z^{-n}$.

Using the  relations \r{gln-com} on $F_i(t)$ we can calculate
the residues in \r{rec-f1} and  obtain the following expressions  for  $\ff_{j,i}(t)$, $i<j$:
\begin{equation}\label{res-in}
\ff_{j,i}(t)=(q-q^{-1})^{j-i-1}\ff_i(t) \ff_{i+1}(t)\cdots \ff_{j-1}(t)\,.
\end{equation}
For example, $F_{i+1,i}(t)=F_i(t)$, and
$F_{i+2,i}(t)=(q-q^{-1})F_i(t)F_{i+1}(t)$. The last product is well-defined
according to the analyticity properties of the product $F_i(t)F_{i+1}(w)$,
described above. In a similar way, one can show inductively that the product
in the right hand side of \r{res-in} makes sense for any $i<j$.
Formulas \r{res-in} prove that the defining relations for the composed currents
\r{rec-f1} or \r{rec-f111}  yields the same answers for all possible
values $i<a<j$.

Calculating formal integrals in \r{rec-f111} we obtain
the following presentation
 for the composed currents:
\begin{equation}\label{rec-f112a}
F_{j,i}(t)=
F_{j,a}(t)F_{a,i}[0]-qF_{a,i}[0]F_{j,a}(t)+
(q-q^{-1})\sum_{k\leq 0}F_{a,i}[k]\,F_{j,a}(t)\,t^{-k}\,,
\end{equation}
\begin{equation}\label{rec-f112b}
F_{j,i}(t)=
F_{j,a}[0]F_{a,i}(t)-q^{-1}F_{a,i}(t)F_{j,a}[0]+
(q-q^{-1})\sum_{k> 0}F_{a,i}(t)\,F_{j,a}[k]\,t^{-k}\,,
\end{equation}
which is useful for the calculation of their projections.

The analytical properties of the products of the
composed currents, used in the paper, are presented in Appendix~\ref{anal-pr1}.

\bigskip

For two sets of variables $\{u_1,...,u_k\}$ and $\{v_1,...,v_k\}$
we introduce the series
\begin{equation}\label{rat-Y}
\begin{split}
\ds Y(u_k,\ldots,u_1;v_k,\ldots,v_1)&=\ds \prod_{m=1}^k\frac{1}{1-v_m/u_m}
\prod_{m'=m+1}^{k}\frac{q-q^{-1}v_{m'}/u_m}{1-v_{m'}/u_m}\\
&=\ds \prod_{m=1}^k\frac{1}{1-v_m/u_m}
\prod_{m'=1}^{m-1}\frac{q-q^{-1}v_m/u_{m'}}{1-v_m/u_{m'}}\ .
\end{split}
\end{equation}

Consider again a collection of segments $\seg{\bar\rr}{\bar\lll}$ and  associated
 set of variables
$\bar t_{\seg{\bar\rr}{\bar\lll}}$.
Let $j=\textrm{max}(a)$ such that $\rr_b=\lll_b$ for $b=a+1,\ldots,N-1$.
Let $\bar\ss$ be a set of nonnegative integers, which besides the inequalities
\begin{equation}\label{rest}
0\leq \ss_a\leq \rr_a-\lll_a\,,\quad a=1,\ldots,j
\end{equation}
satisfies the following admissibility conditions
\begin{equation}\label{non-inc}
0=\ss_0\leq\ss_1\leq\ss_2\ldots\leq \ss_{j-1}\leq\ss_{j}=
\rr_{j}-\lll_j.
\end{equation}

Having the set $\bar\ss$ which satisfy both restrictions
 \r{rest} and \r{non-inc} we define a series depending on the set of the
 variables  $\bar t_{\seg{\bar\rr}{\bar\rr-\bar\ss}}$:
\begin{equation}\label{rat-X}
X(\bar t_{\seg{\bar\rr}{\bar\rr-\bar\ss}})=
\prod_{a=1}^{j-1}
Y(t^{a+1}_{\rr_{a+1}-s_{a+1}+s_{a}},\ldots,t^{a+1}_{\rr_{a+1}-s_{a+1}+1};
t^{a}_{\rr_{a}},\ldots,t^{a}_{\rr_{a}-s_{a}+1} )
\end{equation}
When $j=1$ we set $X(\cdot)=1$.

Define a special ordered product of the composed currents, which we call {\it a string}:
\begin{equation}\label{string}
\F^j_{\bar\ss}({\bar t}^j_{\seg{\rr_j}{\lll_j}})=
\prod^{\longleftarrow}_{j\geq a\geq 1}\sk{
\prod^{\longleftarrow}_{\lll_j+\ss_{a}\geq \ell> \lll_j+\ss_{a-1}}
F_{j+1,a}(t^j_\ell)}\,.
\end{equation}
The string \r{string} depends
 only on the variables $\{ t_{\lll_j+1}^j,\ldots,t^j_{\rr_j}\}$ of the type $j$,
 corresponding to the segment $\seg{\rr_j}{\lll_j}$.
The set of nonnegative integers $\bar\ss$ satisfying the admissibility condition
\r{non-inc} divides the segment $\seg{\rr_j}{\lll_j}$ into $j$ subsegments
$\seg{\lll_j+\ss_a}{\lll_j+\ss_{a-1}}$ for $a=1,\ldots,j$.
This division defines the product of the composed currents in the string \r{string}.

Besides the  product \r{string} which we called the string
we   consider the inverse ordered product of the same composed currents
which we call the {\it inverse string}
\begin{equation}\label{string-inv}
\tF^j_{\bar s}({\bar t}^j_{\seg{\rr_j}{\lll_j}})=
\prod^{\longrightarrow}_{1\leq a\leq j}\sk{
\prod^{\longrightarrow}_{\lll_j+\ss_{a-1}< \ell\leq \lll_j+\ss_{a}}
F_{j+1,a}(t^j_\ell)}\,.
\end{equation}
The inverse string satisfies analytical properties
formulated in Appendix~\ref{anal-pr1} which allows
 to calculate the projection of the direct and
inverse string.

\subsection{Recurrence relation}\label{recur}

Let $\bar n$  be the set of nonnegative integers
$\bar n=\{n_{1},\ldots, n_{N-1}\}$.
We claim that the projection \r{uwf} satisfies the recurrence relation given
by the following
\begin{proposition}\label{prop45}
\begin{equation}\label{main-rec}
\begin{split}
\mathcal{W}^{N-1}&(\bar{t}_{\segg{\bar{n}}})= \sum_{\bar s}
\prod_{a=1}^{N-1}\frac{1}{(s_a-s_{a-1})!(n_a-s_{a})!}\times\\[5mm]
&\times\ \qSym_{\ \bar t_{\segg{\bar n}}}
\sk{Z_{\bar\ss}(\bar t_{\segg{\bar n}})\cdot X(\bar t_{\seg{\bar
n}{\bar n-\bar s}})\cdot {\Pfp}\sk{\F^{N-1}_{\bar s}(\bar
t^{N-1}_{\segg{n_{N-1}}})} \cdot
\mathcal{W}^{N-2}(\bar{t}_{\segg{\bar n-\bar s}})},
\end{split}
\end{equation}
where the sum is taken over all collections
$\bar \ss=\{\ss_{1},\ss_{2},\ldots,\ss_{N-1}\}$, such
that
$0=\ss_0\leq\ss_1\leq\ss_2\ldots\leq \ss_{N-2}\leq\ss_{N-1}=n_{N-1}$
and $0\leq s_a\leq n_a$, $a=1,\ldots,N-2$.
\end{proposition}

Due to the restriction $s_{N-1}=n_{N-1}$ and definition \r{Zser},
the series $Z_{\bar\ss}(\bar t_{\segg{\bar n}})$
depends only on the variables with type $N-1$
missing.
The same restriction
implies that the set $\bar{t}_{\segg{\bar n-\bar s}}$
does not contain variables of the type $N-1$ and an element
$\mathcal{W}^{N-2}(\bar{t}_{\segg{\bar n-\bar s}})$ is a universal weight function for the
algebra $\Uqgl{N-1}$. This fact allows to iterate the
recurrence relation \r{main-rec}  and reduce the universal weight function \r{uwf}
to the $q$-symmetrization of the product of projection of strings. This will be done in the
next section.
\medskip

{\it  Proof of Proposition \ref{prop45}.}\ \ The proof of the
recurrence relation \r{main-rec} follows the strategy of calculation
of the universal weight function which was described in the
introductory part to the Section~\ref{sec4}.

Set $\bar n'=\{n_1,\ldots, n_{N-2},0\}$.
We have a decomposition
$\F(\bar t_{\segg{\bar n}})=
\F(\bar t^{N-1}_{\segg{ n_{N-1}}})\F(\bar t_{\segg{\bar n'}})$
where the first term  only contains the variables of the type $N-1$,
and the last term does not contain the variables of the type $N-1$
\begin{equation*}
\begin{split}
\F(\bar t^{N-1}_{\segg{ n_{N-1}}})&=
\ff_{N-1}(t^{N-1}_{n_{N-1}})\cdots \ff_{N-1}(t^{N-1}_{1}),
\\
\F(\bar t_{\segg{\bar n'}})&=
\ff_{N-2}(t^{N-2}_{n_{N-2}})\cdots \ff_{N-2}(t^{N-2}_{1})\cdots
\ff_{1}(t^{1}_{n_{1}})\cdots
\ff_{1}(t^{1}_{1})\,.
\end{split}
\end{equation*}
We apply  to the product $\F(\bar t_{\segg{\bar n'}})$  the  ordering procedure of
Proposition~\ref{decomp} and substitute the result into \r{uwf}:
\begin{equation}\label{rr1}
\begin{split}
\mathcal{W}^{N-1}(\bar t_{\segg{\bar n}})&=
\sum_{s_{N-2},\ldots,s_1}\prod_{a=1}^{N-2}\frac{1}{s_a!(n_a-s_a)!}
\ \qSym_{\ \bar t_{\segg{\bar n'}}}
\left(Z_{\bar s'}(\bar t_{\segg{\bar n'}})\right.\times\\
&\left.\times\ \Pfp \left( \F(\bar t^{N-1}_{\segg{n_{N-1}}})
 \Pfm \sk{\F(\bar t_{\seg{\bar n'}{\bar n'-\bar s'}})}\right)
\cdot \mathcal{W}^{N-2}(\bar t_{\segg{\bar n'-\bar s'}})\right).
\end{split}
\end{equation}
The sum is taken over all nonnegative integers $\{s_1,\ldots,s_{N-2}\}$
 such that
 $s_a\leq n_a $, $a=1,\ldots,N-2$; $\bar s'$ means the collection
 $\{s_1,\ldots,s_{N-2},0\}$ and
$q$-symmetrization is performed over the variables
$\bar t_{\segg{\bar n'}}$.

\begin{lemma}\label{lemma43a}
For any $\bar n=\{n_{1},\ldots,n_{N-1}\}$, ${\bar
s}=\{s_1,\ldots,s_{N-1}\}$ such that all $s_a\leq n_a$ for
$a=1,\ldots,N-2$ and $s_{N-1}=n_{N-1}$ we have
\begin{equation}\label{rr2a}
\begin{split}
\Pfp&\sk{ \F(\bar t^{N-1}_{\segg{n_{N-1}}})
 \Pfm \sk{\F(\bar t_{\seg{\bar n'}{\bar n'-\bar s'}})}}=\\
 \ds
&=\frac{1}{s_1!}\prod_{a=2}^{N-1}\frac{1}{(s_a-s_{a-1})!}\
\qSym_{\ \bar t_{\seg{\bar n}{\bar n-\bar s}}}
\sk{X(\bar t_{\seg{\bar n}{\bar n-\bar s}})\cdot
\Pfp\sk{\F^{N-1}_{\bar s}(\bar t^{N-1}_{\segg{n_{N-1}}})}}
\end{split}
 \end{equation}
 if $s_{1}\leq s_{2}\leq\cdots\leq s_{N-1}$; otherwise the projection in
  \rf{rr2a} is equal to zero.
\end{lemma}
In \rf{rr2a} we follow the above notations
$\bar n'=\{n_{1},\ldots,n_{N-2},0\}$, ${\bar s'}=\{
s_1,\ldots,s_{N-2},0\}$. The series $X(\bar t_{\seg{\bar n}{\bar n-\bar s}})$ is
defined in \rf{rat-X}.  The proof of Lemma is given in the end of
this section.
\smallskip

Substitute \r{rr2a} into \r{rr1}.  By definition \r{Zser} the series
$Z_{\bar s'}(\bar t_{\segg{\bar n'}})=Z_{\bar s}(\bar t_{\segg{\bar
n}})$ is symmetric with respect to permutations of the variables
$t_{\seg{\bar n}{\bar n-\bar s}}$
of the same type,
 and the universal weight function
$\mathcal{W}^{N-2}(\bar t_{\segg{\bar n'-\bar s'}})=
\mathcal{W}^{N-2}(\bar t_{\segg{\bar n-\bar s}})$ does now depend on
the variables $t_{\seg{\bar n}{\bar n-\bar s}}$.
 We can so include the  series $Z_{\bar s}(\bar t_{\segg{\bar n}})$
and $\mathcal{W}^{N-2}(\bar t_{\segg{\bar n-\bar s}})$ inside the
$q$-symmetrization $\qSym_{\ \bar t_{\seg{\bar n}{\bar
n-\bar s}}}(\cdot)$ and replace the double symmetrization by a
single one, using \r{exa3}:
\begin{equation*}
\qSym_{\ \bar t_{\segg{\bar n'}}}\ \qSym_{\
\bar t_{\seg{\bar n}{\bar n-\bar s}}}(\ \cdot\ )=
\prod_{a=1}^{N-2}s_a!\ \qSym_{\ \bar t_{\segg{\bar n}}}(\ \cdot\ )
\end{equation*}
Proposition~\ref{prop45} is proved.
\hfill$\square$
\medskip

\noindent {\it Proof of Lemma \ref{lemma43a}}. For any $j=1,...,N-1$
denote by $U_{j}$ the subalgebra of $\U_f$ formed by the modes of
$\ff_{1}(t),\ldots,\ff_{j}(t)$. Let $U_j^\varepsilon=U_j\cap
{\mathrm Ker}\, \varepsilon$ be the corresponding augmentation
ideal.

We claim first that the projection $\Pfm \sk{\F(\bar t_{\seg{\bar
n'}{\bar n'-\bar s'}})}$ can be presented as
\begin{equation}\label{rr6}
\begin{split}
 &\overline{\rm Sym}_{\ \bar t_{\seg{\bar n'}{\bar n'-\bar s'}}}
\sk{\frac{1}{s_1!}\frac{X(\bar t_{\seg{\bar n'}{\bar n'-\bar s'}})}
{\prod_{a=1}^{N-3}(s_{a+1}-s_a)!}\
\ \Pfm\sk{\F^{N-2}_{\bar s}(\bar t^{N-2}_{\seg{n_{N-2}}
{n_{N-2}-s_{N-2}}})  }}
\\ &\qquad\qquad\qquad \mod\  \Pfm\sk{ U_{N-3}^\varepsilon}\cdot
U_{N-1}
\end{split}
\end{equation}
 if admissibility conditions
$0\leq s_1\leq s_2\leq\cdots\leq s_{N-3}\leq s_{N-2}$
are satisfied and is zero modulo $\Pfm\sk{ U_{N-3}^\varepsilon}\cdot
U_{N-1}$ otherwise.

This  can be shown by  iteratively using  Proposition \ref{fact3},
 proved in Appendix~\ref{pr6.7}.
 Due to  \r{pgln} and \r{pr-prop} we have
$\Pfm(\Ff_1\cdot \Ff_2)=\Pfm(\Ff_1\cdot \Pfm(\Ff_2))$ for any elements
$\Ff_1,\Ff_2\in\U_f$.
 Thus we can present the projection
 $\Pfm \sk{\F(\bar t_{\seg{\bar n'}{\bar n'-\bar s'}})}$
 as
 \begin{equation}\label{rr44}
 \begin{split}
\frac{1}{s_1!}\Pfm \!\!\left(\prod^{\longleftarrow}_{N-2\geq a\geq3}\!\! \F(\bar
t^{a}_{\seg{n_a}{n_a-s_a}})\cdot
  \F(\bar t^2_{\seg{n_2}{n_2-s_2}})
\cdot \ \qSym_{\ \bar
t^1_{\seg{n_1}{n_1-s_1}}}\!\! \Pfm\!\sk{\F(\bar
t^1_{\seg{n_1}{n_1-s_1}})}\right)\!.
\end{split}
\end{equation}
 We now apply
\r{pr-ap1} with $j=2$ to the last two terms of \rf{rr44}, and under
condition $s_1\leq s_2$ replace them by
\begin{equation}\label{rr45}
\begin{split}
&\frac{1}{(s_2-s_1)!}\ \qSym_{\ {\bar t}_{\seg{\bar n^{(2)}}{{\bar n}^{(2)}-
{\bar s}^{(2)}}}}\left(X({\bar t}_{\seg{{\bar n}^{(2)}}{{\bar n}^{(2)}-
{\bar s}^{(2)}}})\F^2_{{\bar s}^{(2)}}(\bar t^2_{\seg{n_2}{n_2-s_2}})
\right)=\\
&\frac{1}{(s_2-s_1)!}\ \qSym_{\ \bar
t^2_{\seg{n_2}{n_2-s_2}},\bar t^1_{\seg{n_1}{n_1-s_1}}}
\Bigl(Y(t^2_{n_2-s_2+s_1},\ldots,t^2_{n_2-s_2+1};t^1_{n_1},\ldots,t^1_{n_1-s_1+1})\\
&\qquad\times \ff_{32}(t^2_{n_2})\cdots\ff_{32}(t^2_{n_2-s_2+s_1+1})
\ff_{31}(t^2_{n_2-s_2+s_1})\cdots\ff_{31}(t^2_{n_2-s_2+1})\Bigr)\,.
\end{split}
\end{equation}
 modulo
$\Pfm\sk{ U_{1}^\varepsilon}\cdot U_{2}$.
Here we use the notation ${\bar n}^{(2)}=\{n_1,n_2,0,\ldots,0\}$ and
${\bar s}^{(2)}=\{s_1,s_2,0,\ldots,0\}$.
If $s_2<s_1$, the last two terms in  \rf{rr44} are zero modulo
 $\Pfm\sk{ U_{1}^\varepsilon}\cdot U_{2}$.
Due to the commutativity $[\ff_{a}(t),\ff_{1}(t')]=0$ for $a\geq 3$
 we can  move the elements of
 $\Pfm\sk{ U_{1}^\varepsilon}$
 to the left
through the first product $\prod\nolimits_{a\geq 3}\F(\bar
t^{a}_{\seg{n_a}{n_a-s_a}})$ and then out of the projection
$\Pfm\,$, since
 $\Pfm(\Pfm(\Ff')\cdot
\Ff)=\Pfm(\Ff')\cdot \Pfm(\Ff)$. These terms are absorbed into
$\Pfm\sk{ U_{N-3}^\varepsilon}\cdot U_{N-1}$ in \rf{rr6}.

Going further, we replace the appearing string
$\F^2_{{\bar s}^{(2)}}(\bar t^2_{\seg{n_2}{n_2-s_2}}) $   by
its projection $\Pfm$ and apply Proposition \ref{fact3} to the product
$$\F(\bar t^{3}_{\seg{n_3}{n_3-s_3}})\cdot
\qSym_{\ {\bar t}_{\seg{\bar n^{(2)}}{{\bar n}^{(2)}-
{\bar s}^{(2)}}}}\left(X({\bar t}_{\seg{{\bar n}^{(2)}}{{\bar n}^{(2)}-
{\bar s}^{(2)}}})\Pfm\sk{\F^2_{{\bar s}^{(2)}}
(\bar t^2_{\seg{n_2}{n_2-s_2}})}
\right).$$
In  the notation ${\bar n}^{(3)}=\{n_1,n_2,n_3,0,\ldots,0\}$ and
${\bar s}^{(3)}=\{s_1,s_2,s_3,0,\ldots,0\}$, and the assumption $s_2\leq s_3$
this product is replaced  by
$$\frac{1}{(s_3-s_2)!}\ \qSym_{\ {\bar t}_{\seg{\bar n^{(3)}}{{\bar n}^{(3)}-
{\bar s}^{(3)}}}}\left(X({\bar t}_{\seg{{\bar n}^{(3)}}{{\bar n}^{(3)}-
{\bar s}^{(3)}}})\F^3_{{\bar s}^{(3)}}(\bar t^3_{\seg{n_3}{n_3-s_3}})
\right)$$
modulo elements of
$\Pfm\sk{ U_{2}^\varepsilon}\cdot U_{3}$, which are again moved
to the left out of the projection and are absorbed into
$\Pfm\sk{ U_{N-3}^\varepsilon}\cdot U_{N-1}$ in \rf{rr6}.
Finally we get \rf{rr6}.

For the calculation of the projection $\Pfp \sk{
\F(\bar t^{N-1}_{\segg{n_{N-1}}})
 \Pfm \sk{\F(\bar t_{\seg{\bar n'}{\bar n'-\bar s'}})}}$  we replace
 the second factor
 by \r{rr6}. Elements of $\Pfm\sk{ U_{N-3}^\varepsilon}\cdot
U_{N-1}$ do not contribute, since  any element of $U_{N-3}$ commutes with
modes of  $\ff_{N-1}(t)$:
$\Pfp \Big(
\F(\bar t^{N-1}_{\segg{n_{N-1}}})
\cdot \Pfm\sk{ U_{N-3}^\varepsilon}\cdot
U_{N-1} \Big)=0\,.$
Thus we are rest to calculate (we set $s_0\equiv 0$)
\begin{equation*}
\prod_{a=1}^{N-2}\frac{1}{(s_{a}-s_{a-1})!}\
\qSym_{\ \bar t_{\seg{\bar n'}{\bar n'-\bar s'}}}\
\!\Pfp\!\Bigl(
X(\bar t_{\seg{\bar n'}{\bar n'-\bar s'}})
\F(\bar t^{N-1}_{\segg{n_{N-1}}}) \Pfm\!\sk{\F^{N-2}_{\bar s}(\bar t^{N-2}_{\seg{n_{N-2}}
{n_{N-2}-s_{N-2}}})  }\Bigr)\!.
\end{equation*}
Due to Proposition \ref{fact3} the latter expression
is  non-zero only iff $s_{N-2}\leq n_{N-1}$ and is equal to the
right hand side of \r{rr2a}. \hfill{$\square$}

\subsection{Iteration of the recurrence relation}

 Let $\admis{s}=\{s_i^j,\ 1\leq i\leq j\}$ be a triangular matrix with nonnegative integer coefficients.
We say that the matrix $\admis{s}$ is $\bar{n}$-admissible, and denote this
by the symbol $\admis{s}\prec \bar{n}$, if it does not
increase in the lines and its sum over the columns is $\bar{n}$:
\begin{equation}\label{admis}
n_a=\sum_{b=a}^{N-1}s_a^b\,,\qquad a=1,\ldots, N-1\,.
\end{equation}
We also follow the convention $s_0^j=0$ for $j=1,...,N-1$.
\begin{equation}\label{ad-mat}
\admis{s}=
\left(
\begin{array}{ccccc}
s_{1}^{1}&&&&\\[2mm]
s^{2}_1&s^{2}_{2}&&0&\\[2mm]
\vdots&\vdots&\ddots&&\\[2mm]
 s_{1}^{N-2}&s_{2}^{N-2}&\ldots&s_{N-2}^{N-2}&\\[2mm]
s_{1}^{N-1}& s_{2}^{N-1}&\ldots&s_{N-2}^{N-1}&s_{N-1}^{N-1}
\end{array}
\right)\quad
\begin{array}{l}
0=s^1_0\leq s_{1}^{1} \\[2mm]
0=s^2_0\leq s^{2}_{1}\leq s^{2}_2\\[2mm]
\vdots\\[2mm]
0=s_0^{N-2}\!\!\leq s_{1}^{N-2}\leq \ldots \!\leq s_{N-2}^{N-2}   \\[2mm]
0=s_0^{N-1}\!\!\leq s_{1}^{N-1}\leq \ldots\,\leq s_{N-1}^{N-1}
\end{array}
\end{equation}

Let $\bar s^j$, $j=1,\ldots,N-1$ be the $j$-th line of the admissible matrix $\admis{s}$.
Define a collection of vectors
\begin{equation}\label{si-def}
\bar\pp(\bar s)^j=\bar s^{j}+\bar s^{j+1}+\cdots+\bar s^{N-2}+\bar s^{N-1}\,,\quad
 j=1,\ldots,N-1
\end{equation}
with non-negative integer components.
Set $\bar\pp(\bar s)^N=\bar 0$. Note that according to admissibility
condition \r{admis} $\bar\pp(\bar s)^1=\bar n$.

The iteration of the recurrence relations \r{main-rec}
gives the following
\begin{theorem}\label{main-th}
The weight function \r{uwf} can
be presented as a total $q$-sym\-metri\-za\-tion  of the sum over all
 $\bar{n}$ admissible matrices
$\admis{s}$ of the ordered products of the projections of strings with rational coefficients:
\begin{equation}\label{Wt1}
\begin{split}
&\ds \mathcal{W}^{N-1}(\bar{t}_{\segg{\bar{n}}]})=  \qSym
_{\ \bar t_{\segg{\bar n}}}
\sum_{\admis{s}\prec\, \bar{n}}\left(\prod_{b=1}^{N-1}\prod_{a=1}^b
\frac{1}{(s^b_a-s^b_{a-1})!}
\right.\times\\ &\quad \times\ds \left.
\prod_{j=3}^{N-1}
Z_{\bar s^j}(\bar t_{\segg{\bar n-\bar\pp(\bar s)^{j+1}}})
\prod_{j=2}^{N-1}
X(\bar t_{\seg{\bar n-\bar\pp(\bar s)^{j+1}}{\bar n-\bar\pp(\bar s)^{j}}})
\prod_{N-1\geq j\geq 1}^{\longleftarrow}
\Pfp\sk{\F_{\bar s^j}^j (\bar t^j_{\segg{s^j_j}})}
\right).
\end{split}
\end{equation}
\end{theorem}

The rational series  in \r{Wt1}
may be gathered into a single multi-variable series. Set
\begin{equation}\label{gen-ser}
\begin{split}
\mathcal{Z}_{\admis{s}}&(\bar t_{\segg{\bar n}})=
\prod_{j=3}^{N-1}
Z_{\bar s^j}(\bar t_{\segg{\bar n-\bar\pp(\bar s)^{j+1}}})
\prod_{j=2}^{N-1}
X(\bar t_{\seg{\bar n-\bar\pp(\bar s)^{j+1}}{\bar n-\bar\pp(\bar s)^{j}}})=\\
&=\prod_{b=2}^{N-1}\prod_{a=1}^{b-1}\prod_{\ell=1}^{s^b_a}
\frac{1}{1-t^{a}_{\ell+n_a-\pp^{b}_{a}}/t^{a+1}_{\ell+n_{a+1}-\pp^{b}_{a+1}}}
\prod_{\ell'=1}^{\ell+n_{a+1}-\pp^b_{a+1}-1}
\frac{q-q^{-1}t^{a}_{\ell+n_a-\pp^{b}_{a}}/t^{a+1}_{\ell'}}
{1-t^{a}_{\ell+n_a-\pp^{b}_{a}}/t^{a+1}_{\ell'}}
\end{split}
\end{equation}
where $\pp^b_a$ are components of the vector $\bar\pp(\bar s)^b$.
Using \r{gen-ser} the formula for the universal weight
 function \r{Wt1} can be written
in a compact form:
\begin{equation}\label{Wt-fin}
\ds \mathcal{W}^{N-1}(\bar{t}_{\segg{\bar{n}}})=  \qSym
_{\ \bar t_{\segg{\bar n}}}
\sum_{\admis{s}\prec\, \bar{n}}\left(
\frac{\mathcal{Z}_{\admis{s}}(\bar t_{\segg{\bar n}})
}{\prod_{a\leq b}(s^b_a-s^b_{a-1})!}
\prod^{\longleftarrow}_{N-1\geq j\geq 1}
\Pfp\sk{\F_{\bar s^j}^j (\bar t^j_{\segg{s^j_j}})    }
\right).
\end{equation}
 Theorem~\ref{main-th} reduces the calculation of the universal weight function
to the calculation of the projections of the strings.

\subsection{Projection of composed currents and strings}

Let $\Delta$ be the standard comultiplication in $\Uqgln$ defined by \r{coprL}. For any
elements $x,y\in\Uqgln$ we define the adjoint action by the relation
\begin{equation}\label{ad-stan}
{\rm ad}_x\ y= \sum\nolimits_l a(x'_l)\cdot y\cdot x''_l\ , \quad\mbox{where}\quad
\Delta x =\sum\nolimits_l x'_l\ot x''_l
\end{equation}
and $a(x)$ is an antipode map related to the comultiplication \r{coprL}.
By  {\it screening operators} we  understand the operators of the adjoint actions
of  zero modes $F_i[0]$ of  $F_{i}(t)$:
\begin{equation}\label{ad-stan1}
S_i\, (y)={\rm ad}_{F_{i}[0]}\ (y)=
y\ F_{i}[0] - F_{i}[0]\ k^{-1}_{i+1} k_i\ y\  k^{-1}_i k_{i+1}
\end{equation}
where $k_i$ are zero modes of $k_i^+(t)$:  $k_i=k_i^+[0]=k_i^-[0]^{-1}$.

\begin{proposition}\label{com-cur2}For any  $i,j$, $1\leq i<j\leq N $
we have the equalities
\begin{equation}\label{cc-pr}
\begin{split}
\Pfp\sk{\ff_{j+1,i}(u)}&=S_{i}
                       S_{{i+1}}
                       \cdots
                       S_{{j-1}} \sk{ \Pfp( {F}_{j}(u)) }\\
                       &=
                       S_{i}
                       S_{{i+1}}
                       \cdots
                       S_{{j-1}}\left(   {F}_{j}(u) ^{(+)}\right)=
                       \left(S_{i}
                       S_{{i+1}}
                       \cdots
                       S_{{j-1}} \sk{  {F}_{j}(u) }\right)^{(+)}
                       \end{split}
\end{equation}
\end{proposition}
\noindent
{\it Proof.}\  Taking \r{rec-f112a}
for $a=i+1$ and using the definition of the projection $\Pfp$ we obtain
\begin{equation}\label{pro-in}
\Pfp\sk{\ff_{j,i}(t)}\,=\,
S_{i}\bigl(P^+(\ff_{j,i+1}(t))\bigr)\,,\qquad i<j-1\,.
\end{equation}
The first equality  of Proposition follows from \r{pro-in}
by induction. One should take into account
the commutativity of the projections and screening operators proved
in \cite{KP}. Then we apply \rf{DFinverse}.
\hfill$\square$

An analog of Proposition \ref{com-cur2} for $\Pfm\sk{\ff_{j+1,i}(u)}$ is given in
Appendix B.

\medskip

For a set $\{u_1,\ldots,u_n\}$ of  formal variables we introduce a
set of the rational functions
\begin{equation}\label{rat-f}
\varphi_{u_m}(u;u_1,\ldots,u_n)=\prod_{k=1,\ k\neq m}^{n}
\frac{u-u_k}{u_m-u_k}\prod_{k=1}^{n}\frac{q^{-1}u_m-qu_k}{q^{-1}u-qu_k}
\end{equation}
satisfying the normalization conditions
$\varphi_{u_m}(u_s;u_1,\ldots,u_n)=\delta_{ms}$. We set
\begin{equation}\label{long-cur}
\ff_{j+1,i}(u;u_{1},\ldots,u_n)=\ff_{j+1,i}(u)-
\sum_{m=1}^n \varphi_{u_m}(u;u_{1},\ldots,u_n)\ff_{j+1,i}(u_m)
\end{equation}
for $1\leq i\leq j<N$. By \rf{cc-pr} we have
\begin{equation}\label{cc-pr1}
\begin{split}
\Pfp\sk{\ff_{i,j+1}(u;u_{1},\ldots,u_n)}&=S_{i}
                       S_{i+1}
                       \cdots
                       S_{j-1}\sk{F_{j}(u)^{(+)}}-\\
-&\sum_{m=1}^n \varphi_{u_m}(u;u_{1},\ldots,u_n)S_i
S_{i+1}\cdots                       S_{j-1}\sk{F_{j}(u_m)^{(+)}}\ .
\end{split}
\end{equation}
Proposition~\ref{prop1} below and the relation
\r{cc-pr1} suggest a factorized form for the projection of the
inverse string $\tF_{\bar s}(\bar t^j_{\seg{\rr_j,\lll_j}{}})$.
\begin{proposition}\label{prop1}
\begin{equation}\label{pr-st-inv}
\Pfp\sk{\tF^j_{\bar s}(\bar t^j_{\seg{\rr_j}{\lll_j}})}=
\prod^{\longrightarrow}_{1\leq a\leq j}\sk{
\prod^{\longrightarrow}_{\lll_j+s_{a-1}< \ell\leq \lll_j+s_{a}}
\Pfp\sk{F_{j+1,a}(t^j_\ell;t^j_{\lll_j+1},\ldots,t^j_{\ell-1})}}.
\end{equation}
\end{proposition}

\noindent
{\it Proof}\ \ of  Proposition  \rf{pr-st-inv} is shifted to the Appendix~\ref{pr-st3}.
\medskip

Projections of the string \r{string} and of the inverse string \r{string-inv}
are related, namely
\begin{equation}\label{pr-st1}
\begin{split}
&\ds \Pfp\sk{\F^j_{\bar s}(\bar t^j_{\seg{\rr_j}{\lll_j}})}=
\Pfp\sk{\tF^j_{\bar s}(\bar t^j_{\seg{\rr_j}{\lll_j}})}\ \times\\
\\ &\ds\qquad\times
\prod_{\lll_j<\ell<\ell'\leq \rr_j}
\frac{q-q^{-1}t^j_{\ell'}/t^j_{\ell}}{1-t^j_{\ell'}/t^j_{\ell}}
\prod_{1\leq a\leq j}\sk{
\prod_{\lll_j+s_{a-1}<\ell<\ell'\leq \lll_j+s_a}
\frac{1-t^j_{\ell'}/t^j_{\ell}}{q^{-1}-qt^j_{\ell'}/t^j_{\ell}} }.
\end{split}
\end{equation}   This statement is the direct consequence of
 the  relations between composed
currents given by Proposition~\ref{prop4.2}. For the details see
 \cite{KP}.

The particular case of basic relations \rf{gln-com},
$$
k^+_{j}(t')F_{j,j-1}(t)k^+_{j}(t')^{(-1)}=\frac{q-q^{-1}t/t'}{1-t/t'}
F_{j,j-1}(t)
$$
implies the following  relation for the projections:
\begin{equation}\label{gc1}
k^+_{j}(t')\Pfp\sk{F_{j,j-1}(t)}k^+_{j}(t')^{-1}=\frac{q-q^{-1}t/t'}{1-t/t'}
\Pfp\sk{F_{j,j-1}(t;t')}\,,
\end{equation}
and in general
\begin{equation}\label{gc2}
\prod_{a=1}^n k^+_{j}(t_a) \Pfp\sk{F_{j,j-1}(t)}
\prod_{a=1}^n
k^+_{j}(t_a)^{-1}=\prod_{a=1}^n\frac{q-q^{-1}t/t_a}{1-t/t_a}\
\Pfp\sk{F_{j,j-1}(t;t_1,\ldots,t_n)}.
\end{equation}
Applying  a sequence of the screening
operators $S_{i}\cdots S_{j-2}$ to  \r{gc2} we obtain
\begin{equation}\label{gc4}
\prod_{a=1}^n k^+_{j}(t_a)\cdot \Pfp\sk{F_{j,i}(t)}\cdot
\prod_{a=1}^n
k^+_{j}(t_a)^{-1}=\prod_{a=1}^n\frac{q-q^{-1}t/t_a}{1-t/t_a}
\Pfp\sk{F_{j,i}(t;t_1,\ldots,t_n)}
\end{equation}
due to the commutativity of any of zero modes
$F_{i}[0],\ldots, F_{j-1}[0]$ with
$k^+_{j+1}(t)$. We rewrite the  relation \r{gc4} in the
form
\begin{equation}\label{gc5}
\Pfp\sk{F_{j,i}(t;t_1,\ldots,t_n)}\cdot \prod_{a=1}^n
k^+_{j}(t_a)= \prod_{a=1}^n\frac{1-t/t_a}{q-q^{-1}t/t_a}
\prod_{a=1}^n k^+_{j}(t_a)\cdot \Pfp\sk{F_{j,i}(t)}.
\end{equation}
Using \rf{gc5} and \rf{pr-st1} we rewrite the relation  \r{pr-st-inv}
and its analog for the string \rf{string} in the following form:
\begin{equation}\label{gc7}
\begin{split}
&\ds \Pfp\sk{\tF^j_{\bar s}(\bar t^j_{\seg{\rr_j}{\lll_j}})}
\cdot
\prod_{\ell=\lll_j+1}^{\rr_j} k^+_{j+1}(t^j_\ell)=\\
 &\ds\quad =
\prod_{\lll_j<\ell<\ell'\leq
\rr_j}\frac{1-t^j_{\ell'}/t^j_\ell}{q-q^{-1}t^j_{\ell'}/t^j_\ell}
\prod^{\longrightarrow}_{1\leq a\leq j}\sk{
\prod^{\longrightarrow}_{\lll_j+s_{a-1}< \ell\leq \lll_j+s_{a}}
\Pfp\sk{F_{j+1,a}(t^j_\ell)}k^+_{j+1}(t^j_\ell)}\,,
\end{split}
\end{equation}
\begin{equation}\label{gc8}
\begin{split}
&\ds \Pfp\sk{\F^j_{\bar s}(\bar t^j_{\seg{\rr_j}{\lll_j}})}
 \prod_{\ell=\lll_j+1}^{\rr_j} k^+_{j+1}(t^j_\ell)
=  \prod_{1\leq a\leq j}\sk{
\prod_{\lll_j+s_{a-1}<\ell<\ell'\leq \lll_j+s_a}
\frac{1-t^j_{\ell'}/t^j_{\ell}}{q^{-1}-qt^j_{\ell'}/t^j_{\ell}}
}\times\\
&\ds\quad\qquad\times \prod^{\longrightarrow}_{1\leq a\leq j}\sk{
\prod^{\longrightarrow}_{\lll_j+s_{a-1}< \ell\leq \lll_j+s_{a}}
\Pfp\sk{F_{j+1,a}(t^j_\ell)}k^+_{j+1}(t^j_\ell)}.
\end{split}
\end{equation}

\section{Weight function and $\LL$-operators}
\subsection{From Gauss coordinates to $\LL$-operator's entries}
The results of the previous section show that the projection of the
string multiplied by certain number of  $k^+_{j+1}(t^j_\ell)$ can be
factorized in such a way that each factor is a product of projection
of the composed currents and of some $k^+_{j+1}(t^j_\ell)$. In this
section we use this observation and express the weight functions via
matrix elements of $\LL$-operators \r{L-op}.
 First, we relate projection of composed currents
with Gauss coordinates of $\LL$-operators. This is given by the following
\begin{proposition}\label{idenGC}
We have for any $ i<j-1$
\begin{equation}\label{ide-pr}
\Pfp\sk{F_{j,i}(t)}\,=\,(q-q^{-1})^{j-i-1}\FF^{+}_{j,i}(t)\,.
\end{equation}
\end{proposition}
\noindent{\it Proof.}\ \,
 We use the equality
\begin{equation}\label{cor1}
(q-q^{-1})\FF^{+}_{j,i}(t)\,=\,S_{i}\sk{\FF^{+}_{j,i+1}(t)}\,,\qquad i<j-1\,.
\end{equation}
It is proved in \cite{KPT} and is a direct consequence of the
 relations \r{L-op-com} taken for the modes of $\LL$-operators.
The claim of Proposition now  follows by induction from
Proposition~\ref{com-cur2} and relation \rf{DFinverse}.
\hfill$\square$
\medskip

 Let $V$ be $\Uqgln$-module with a weight singular vector $v$.
  The relation \r{gc8} and Proposition~\ref{idenGC}
allows to rewrite the corollary \r{Wt-fin} to  Theorem~\ref{main-th}
in the following form:
\begin{equation}\label{Wt2}
\begin{split}
\ds {\bf w}_V^{N-1}(\bar{t}_{\segg{\bar{n}}})&=\ds
\qsym(\bar t_{\segg{\bar n}})\
\qSym
_{\ \bar t_{\segg{\bar n}}}
\sum_{\admis{s}\prec\, \bar{n}}\left(
\frac{(q-q^{-1})^{\sum_{b=1}^{N-1}(n_b-s_b^b)}}
{ \prod_{a\leq b}^{N-1}(s^b_a-s^b_{a-1})!}
\ \mathcal{Z}_{\admis{s}}(\bar t_{\segg{\bar n}})\!\!\!\!\!
\prod_{N-1\geq b\geq 1}^{\longleftarrow}\left(
\prod_{1\leq a\leq b}^{\longrightarrow}
  \right.\right.\\
 &\left.\left.
\times \prod_{s^b_{a-1}<\ell \leq s^b_a}^{\longrightarrow}
\left( \FF^{+}_{b+1,a}(t^b_\ell)k^+_{b+1}(t^b_\ell)
\prod_{\ell'=\ell+1}^{s^b_a}
\frac{t^b_\ell-t^b_{\ell'}}{q^{-1}t^b_\ell-qt^b_{\ell'}} \right)\right)
 \prod_{\ell=s^b_b+1}^{n_b}k^+_{b+1}(t^b_\ell)\right)v
\end{split}
\end{equation}
where the series $\mathcal{Z}_{\admis{s}}(\bar t_{\segg{\bar n}})$ is given by \r{gen-ser}.
\medskip

{}For any $c=1,\ldots,N$ denote by $I_c$ the left ideal of
$U_q(\widehat{\mathfrak{b}}^+)$, generated by the modes of
$\EE^+_{j,i}(u)$ with $ i>j\geq c$. We have inclusions $0=I_N\subset
 I_{N-1}\subset\cdots\subset I_1$.

\begin{lemma}\label{lemmaideal} Fix any $c=1,\ldots,N-1$. Then

\noindent (i) the left ideal $I_c$ is generated by  modes of
$\,\LL^+_{i,j}(u)$ with $ i>j\geq c$;

\noindent (ii)  for any $a$ and $b$ with $a<b$ and $b\geq c$ we have
equalities
\begin{equation}\label{comparison}
\LL^+_{a,b}(t)\,\equiv\, \FF^+_{b,a}(t)k^+_b(t)\, \mod \ I_c,\
\qquad \LL^+_{b,b}(t)\,\equiv\, k^+_b(t)\, \mod \
I_c,
\end{equation}

\noindent (iii) for any  $a\leq c$  and $b\geq c$ the modes of
$\ \LL^+_{a,c}(t)$ and of $\ \LL^+_{b,b}(t)$ normalize the ideal $I_c$:
\begin{equation} \label{normalize}
 I_c\cdot  \LL^+_{a,c}(t) \subset I_c\,\qquad
I_c\cdot  \LL^+_{b,b}(t) \subset I_c\,.
\end{equation}
\end{lemma}

\noindent
{\it Proof.}\  We have three types of relations  \r{ide-pr}:
 \begin{align}\label{GF}
\LL^{+}_{a,b}(t)&=\FF^{+}_{b,a}(t)k^+_{b}(t)+\sum_{b<m\leq N}
\FF^{+}_{m,a}(t)k^+_{m}(t)\EE^{+}_{b,m}(t),\qquad a<b,\\
\label{GK}
\LL^{+}_{b,b}(t)&=k^+_{b}(t)\quad +\qquad\sum_{b<m\leq N} \FF^{+}_{m,b}(t)k^+_{m}(t)
\EE^{+}_{b,m}(t),\\
\label{GE}
\LL^{+}_{a,b}(t)&=k^+_{a}(t)\EE^{+}_{b,a}(t)+\sum_{a<m\leq N}
\FF^{+}_{m,a}(t)k^+_{m}(t)\EE^{+}_{b,m}(t),\qquad   a>b\,.
\end{align}
Denote $\bar{\EE}^+_{j,i}(u)=k^+_{i}(t)\EE^{+}_{j,i}(t)$ and
$\bar{\FF}^+_{i,j}(u)=\FF^{+}_{i,j}(t)k^+_{i}(t)$ for $i>j$. In
these notations the relation \rf{GE} looks as
\begin{equation}\label{GE1}
\LL^{+}_{a,b}(t)=\bar{\EE}^{+}_{b,a}(t)+\sum_{a<m\leq N}
\FF^{+}_{m,a}(t)\bar{\EE}^{+}_{b,m}(t),\qquad   a>b\,.
\end{equation}
 Since $k^+_b(t)$ is
invertible, the ideal $I_c$ is generated by the modes of
$\bar{\EE}^+_{j,i}(u)$ with $ i>j\geq c$  as well. Now we inverse
the relations \rf{GE1}, that is we write $\bar{\EE}^{+}_{j,N}(t)
=\LL^+_{N,j}(t)$ for $j<N$;  then $\bar{\EE}^{+}_{j,N-1}(t)=
\LL^+_{N-1,j}(t)- \FF^{+}_{N,N-1}(t)\LL^+_{N,j}(t)$ for $j<N-1$ and
so on by induction. This proves (i). In the same manner we  rewrite
first \rf{GF} and \rf{GK} in $\bar{\FF}^+_{i,j}(u)$,
${\EE}^+_{j,i}(u)$ and $k^+_{i}(t)$ and prove by induction that
$\LL^+_{a,b}(t)\,\equiv\, \bar{\FF}^+_{b,a}(t)$\ $ \mod \ I_c\ $
{and} $\LL^+_{b,b}(t)\,\equiv\, k^+_b(t)$\ mod \ $I_c$ for $b\geq c$
and $a<b$, which means (ii).

The statement (iii) is a corollary of the Yang-Baxter relations
\rf{YBeq}. Let $\RR_{ij;kl}(u,v)$ be the matrix elements of the
$R$-matrix \rf{UqglN-R}. We have for $a<c<k$ and for $a<c<j<i$:
\begin{align*}
\LL^+_{k,c}(u)\LL^+_{a,c}(t)=
\frac{\RR_{cc;cc}(u,t)}{\RR_{kk;aa}(u,t)}\LL^+_{a,c}(t)\LL^+_{k,c}(u)&-
\frac{\RR_{ka;ak}(u,t)}{\RR_{kk;aa}(u,t)}\LL^+_{a,c}(u)\LL^+_{k,c}(t),
\\
\LL^+_{i,j}(u)\LL^+_{a,c}(t)=
\frac{\RR_{jj;cc}(u,t)}{\RR_{ii;aa}(u,t)}\LL^+_{a,c}(t)\LL^+_{i,j}(u)&-
\frac{\RR_{ia;ai}(u,t)}{\RR_{ii;aa}(u,t)}\LL^+_{a,c}(u)\LL^+_{i,c}(t)
\\&+
\frac{\RR_{cj;jc}(u,t)}{\RR_{ii;aa}(u,t)}\LL^+_{a,c}(t)\LL^+_{i,c}(u).
\end{align*}
These relations precisely mean the inclusion $I_c\cdot
\LL^+_{a,c}(t) \subset I_c$. The second part of (iii) is proved in
an analogous manner and is actually well known.
\hfill$\square$

\begin{theorem}\label{th-Wt3} For any $\Uqgln$ module $V$ with a weight singular
vector $v$ we have
\begin{equation}\label{Wt3}
\begin{split}
&\ds {\bf w}_V^{N-1}(\bar{t}_{\segg{\bar{n}}})=\ds
\qsym(\bar t_{\segg{\bar n}})\
 \qSym_{\ \bar t_{\segg{\bar n}}}
\sum_{\admis{s}\prec\, \bar{n}}\left(
\frac{(q-q^{-1})^{\sum_{b=1}^{N-1}(n_b-s_b^b)}}{\prod_{a\leq b}^{N-1}
(s^b_a-s^b_{a-1})!}
\ \mathcal{Z}_{\admis{s}}(\bar t_{\segg{\bar n}})\times
\right.\\
 &\ds \prod_{N-1\geq b\geq 1}^{\longleftarrow}\left.\!\!\left(
\prod_{1\leq a\leq b}^{\longrightarrow}
\prod_{\ell=s^b_{a-1}+1}^{s^b_a}\sk{\LL^{+}_{a,b+1}(t^b_\ell)
\prod_{\ell'=\ell+1}^{s^b_a}
\frac{t^b_\ell-t^b_{\ell'}}{q^{-1}t^b_\ell-qt^b_{\ell'}} }\right)
\prod_{\ell=s^b_b+1}^{n_b}\LL^{+}_{b+1,b+1}(t^b_\ell)
\right)v
\end{split}
\end{equation}
\end{theorem}
Observe that non-commutative products over $b$ and $a$ run in the opposite directions
and the ordering in the product over $\ell$ is not important now, because of commutativity
of the matrix elements of $\LL$-operators with the same matrix indices.
\smallskip

\noindent
 {\it Proof}. The theorem states that in \rf{Wt2} we can replace
 each entry of
 $\bar{\FF}^+_{b,a}(t)= \FF^{+}_{b,a}(t)k^+_{b}(t)$ by
 $\LL^{+}_{a,b}(t)$ and each entry of $k_b^+(t)$ in the last product
 of \rf{Wt2} by $\LL^{+}_{b,b}(t)$. This is done with a help of Lemma
 \ref{lemmaideal}. Indeed, we can present the right hand side of
 \rf{Wt2} as a linear combination of terms
 $$A_i^NA_i^{N-1}\cdots A_i^{2} B_i v,$$
 where each $A_i^c$ is an ordered product of some
 $\bar{\FF}^+_{c,a}(t_j^{c-1})$, and $B_i$ is a product of some
 $k_b(t_j^{b-1})$. We start from $A_i^N$. Here we have
 $\bar{\FF}^+_{N,a}(t)=\LL^{+}_{a,N}(t)$ by \rf{GF}. By the statement
 (ii) of Lemma \ref{lemmaideal}
each multiplier of $A_i^{N-1}$ can be written as
$\LL^{+}_{a,N-1}(t_j^{N-2})+x_j$ for some $a<N-1$, parameter
$t_j^{N-2}$ and $x_j\in I_{N-1}$. Due to the part (iii) of Lemma
 \ref{lemmaideal}, we can rewrite the whole product in $A_i^{N-1}$
 as $\prod_j \LL^{+}_{a,N-1}(t_j^{N-2}) +y_{N-1}$ with a single
 $y_{N-1}\in I_{N-1}$. Moving further, we replace the product
$A_i^NA_i^{N-1}\cdots A_i^{2} B_i$ by
\begin{equation}\label{tildeprod}
\left(\bar{A}_i^{N}\bar{A}_i^{N-1}\cdots\bar{A}_i^{2}+y_2\right)B_i
v \,,\end{equation}
 where each $\bar{A}_i^{c}$ equals to ${A}_i^{c}$
with all $\bar{\FF}^+_{c,a}(t)$ replaced by $\LL^+_{a,c}(t)$, and
$y_2\in I_2\subset I_1$. We apply finally the last part of Lemma
 \ref{lemmaideal}, (iii) and replace \rf{tildeprod} by the product
$\bar{A}_i^{N}\bar{A}_i^{N-1}\cdots\bar{A}_i^{2}\bar{B}_i v$, where
$\bar{B}_i$ consists of $\LL^+_{b,b}(t_k^{b-1})$ only. This proves
the theorem.\hfill$\square$
\medskip

We can slightly simplify  the formula \r{Wt3} by a
renormalization of  $q$-sym\-me\-tri\-zation (see \cite{TV3}). We set
for a series or function  $G(\bar t_{\segg{\bar n}})$:
\begin{equation}\label{qs1r2}
\tSym_{\ \bar t_{\segg{\bar n}}}
\ G(\bar t_{\segg{\bar n}})= \qsym(\bar t_{\segg{\bar n}})\
\qSym_{\ \bar t_{\segg{\bar n}}} \sk{ G(\bar t_{\segg{\bar n}})}\,.
\end{equation}

 Let $G^{\rm sym}(u_1,\ldots,u_n)$ be any symmetric function of $n$
 variables $u_k$, that is  $G^{\rm sym}(^\si\bar u)=G^{\rm sym}(\bar u)$
for any element $\si$ from the symmetric group $S_n$. Let
$\qsym(\bar u)=\prod_{k<k'}
\frac{q-q^{-1}u_k/u_{k'}}{1-u_k/u_{k'}}$. One can check the
following property of renormalized $q$-symmetrization:
\begin{equation}\label{po-sim}
\frac{1}{n!}\ \ \tSym_{\ \bar u}\sk{\qsym(\bar u)^{-1}G^{\rm sym}(\bar u)}
=\frac{1}{[n]_q!}\ \  \tSym_{\ \bar u}\sk{G^{\rm sym}(\bar u)}\,,
\end{equation}
where $[n]_q=\frac{q^n-q^{-n}}{q-q^{-1}}$,
 and $[n]_q!=[n]_q[n-1]_q\cdots [2]_q[1]_q$.

We can apply the relation \rf{po-sim} to the right hand side of
\r{Wt3} because  in this formula inside the total $q$-symmetrization
there is a symmetric series in the sets of variables
$\{t^{b}_\ell\}$ for $s^b_{a-1}+1\leq \ell\leq s^b_a$,
$a=1,\ldots,b$ and $b=1,\ldots,N-1$. This follows from the
commutativity of matrix elements of $\LL$-operators
$[\LL^+_{a,b}(t),\LL^+_{a,b}(t')]=0$, and from the explicit form of the
series \r{gen-ser}. Restoring this series and denoting
$\ppp^b_a=n_a-\pp^b_a=s^{a}_a+\cdots+s^{b-1}_a$ we formulate the
following corollary of Theorem~\ref{th-Wt3}
\begin{corollary}\label{cor53}
The off-shell Bethe vectors for quantum affine algebra $\Uqgln$ can be written as
\begin{equation}\label{Wt4}
\begin{split}
\ds {\bf w}_V^{N-1}(\bar{t}_{\segg{\bar{n}}})&=\ds  \tSym
_{\ \bar t_{\segg{\bar n}}}
\sum_{\admis{s}\prec\, \bar{n}}\left((q-q^{-1})^{\sum_{b=1}^{N-1}(n_b-s_b^b)}
\prod_{a\leq b}^{N-1}
\frac{1}{[s^b_a-s^b_{a-1}]_q!}\  \ \times \right.\\
&\times\
\prod_{b=2}^{N-1}\prod_{a=1}^{b-1}\prod_{\ell=1}^{s^b_a}
\frac{1}{1-t^{a}_{\ell+\ppp^{b}_{a}}/t^{a+1}_{\ell+\ppp^{b}_{a+1}}}
\prod_{\ell'=1}^{\ell+\ppp^b_{a+1}-1}
\frac{q-q^{-1}t^{a}_{\ell+\ppp^{b}_{a}}/t^{a+1}_{\ell'}}
{1-t^{a}_{\ell+\ppp^{b}_{a}}/t^{a+1}_{\ell'}}
\\
&\times\ds \left.
\prod_{N-1\geq b\geq 1}^{\longleftarrow} \sk{\prod_{1\leq a\leq b}^{\longrightarrow}
\sk{\prod_{\ell=s^b_{a-1}+1}^{s^b_a} \LL^+_{a,b+1}(t^b_\ell) }
\prod_{\ell=s^b_b+1}^{n_b}\LL^+_{b+1,b+1}(t^b_\ell)}
\right)v.
\end{split}
\end{equation}
\end{corollary}

Expression \r{Wt4} for the off-shell Bethe vectors expressed in terms of matrix elements
of $\LL$-operators can be written as the recurrence relation
\begin{equation}\label{rec-bv2}
\begin{split}
{\bf w}_V^{N-1}(\bar{t}_{\segg{\bar{n}}})&=
\sum_{\bar s^{N-1}}\frac{1}{[s^{N-1}_1]_q!}
\prod_{a=1}^{N-2}\frac{(q-q^{-1})^{s_a^{N-1}}}{[s^{N-1}_{a+1}-s^{N-1}_{a}]_q![n_a-s^{N-1}_{a}]_q!}\\
&\ds\quad \times \tSym_{\ \bar t_{\bar n}}
\left(X(\bar t_{\seg{\bar n}{\bar n-\bar s^{N-1}}})\ Z_{{\bar s}^{N-1}}(\bar t_{\segg{\bar n}})
\prod_{a=1}^{N-2}
\prod_{\ell\,=\,n_{a}-s^{N-1}_{a}+1}^{n_{a}}
\lambda_{a+1}(t^{a}_\ell),
\right.\\
&\ds\quad \left.\times \prod^{\longrightarrow}_{1\leq a\leq N-1}
\sk{\prod_{\ell=s^{N-1}_{a-1}+1}^{s^{N-1}_a}
\LL^+_{a,N}(t^{N-1}_\ell) }
\
{\bf w}_V^{N-2}(\bar{t}_{\segg{\bar{n}-\bar s^{N-1}}})
\right)
\end{split}
\end{equation}
where summation in \r{rec-bv2} runs over first row of the admissible matrix $\admis{s}$.
The vector valued function ${\bf w}_V^{N-2}(\bar{t}_{\segg{\bar{n}-\bar s^{N-1}}})$ is
the off-shell Bethe vector for the algebra $\Uqgl{N-1}$ embedded into $\Uqgln$. This embedding
$\phi:\Uqgl{N-1}\hookrightarrow\Uqgln$ can be described on the level of Gauss coordinates
of the corresponding $\LL$-operators:
\begin{equation}\label{emb1}
\begin{split}
\phi\sk{\FF^\pm_{b,a}(t)^{\<N-1\>}}&= \FF^\pm_{b,a}(t)^{\<N\>}\,,\quad
\phi\sk{\EE^\pm_{a,b}(t)^{\<N-1\>}}= \EE^\pm_{a,b}(t)^{\<N\>}\,,\\
\phi\sk{k^\pm_{a}(t)^{\<N-1\>}}&= k^\pm_{a}(t)^{\<N\>}\,,
\quad 1\leq a<b\leq N-1\,.
\end{split}
\end{equation}
Here $\FF^\pm_{b,a}(t)^{\<N-1\>}$ etc., denote Gauss coordinates of
the source $\Uqgl{N-1}$ while $\FF^\pm_{b,a}(t)^{\<N\>}$ etc., the
Gauss coordinate of the target $\Uqgl{N}$.

\subsection{ Evaluation homomorphism}

The quantum affine algebra $\Uqgln$ contains a quantum group
$U_q(\mathfrak{gl}_N)$. It is generated by the zero modes of
$\LL$-operators \r{L-op}.
Using them, we introduce Cartan-Weyl generators
$\ele_{ab}$, $1\leq a,b\leq N$ of $U_q(\mathfrak{gl}_N)$.
 We set $\LL^\pm\equiv \LL^\pm[0]$, and
\begin{equation}\label{LL+}
\LL^+={\mathrm K}\LL^+_0 =
\sk{ \begin{array}{cccc} \ele_{11} &&&\\
&\ele_{22}&0&\\ & 0&\ddots& \\ &&&    \ele_{NN}
\end{array} }
\sk{ \begin{array}{cccc}
1&&&\\ \hba \ele_{12}&1&0&\\ \vdots & \ddots &\ddots&\\
\hba \ele_{1N}&\cdots&\hba \ele_{N-1,N}&1
\end{array} },
\end{equation}
\begin{equation}\label{LL-}
\LL^-\! =\LL^-_0{\mathrm K}^{-1}\!=
\sk{ \begin{array}{cccc}
1&\!\!-{\hba} \ele_{21}&\cdots&\!\!-{\hba} \ele_{N1}\\
&\ddots&\ddots&\vdots \\
&&1&\!\!-{\hba} \ele_{N,N-1}\!\\
&&&1\\
\end{array} }
\sk{ \begin{array}{cccc}\! \ele^{-1}_{11} &&&\\
&\ele^{-1}_{22}&0&\\ & 0&\ddots& \\   &&&  \!\!  \ele^{-1}_{NN}\!
\end{array} },
\end{equation}
where $\hba=q-q^{-1}$.
Generators $\ele_{aa}$, $\ele_{a,a+1}$ and $\ele_{a+1,a}$ may be considered as
Chevalley generators of $\Uqqq$ with commutation relations
\begin{equation}\label{uq-com}
\ele_{aa}\ele_{bc}\ele^{-1}_{aa}=q^{\d_{ab}-\d_{ac}}\ele_{bc}\,,\qquad
[\ele_{a,a+1},\ele_{b+1,b}]=\d_{ab}\frac{\ele_{aa}\ele^{-1}_{a+1,a+1}-
\ele^{-1}_{aa}\ele_{a+1,a+1}}{q-q^{-1}},
\end{equation}
\begin{equation}
\begin{split}\label{serregl}
\ele^2_{a\pm 1,a}\ele_{a,a\mp 1}-(q+q^{-1})\ele_{a\pm
1,a}\ele_{a,a\mp 1}\ele_{a\pm 1,a}
+\ele_{a,a-1}\ele^2_{a\pm 1,a}&=0,\\
\ele^2_{a,a\pm 1}\ele_{a\mp 1,a}-(q+q^{-1})\ele_{a,a\pm 1}\ele_{a\mp
1,a}\ele_{a,a\pm 1}
+\ele_{a\mp 1,a}\ele^2_{a,a\pm 1}&=0.
\end{split}
\end{equation}

The rest of $\ele_{ab}$ may be constructed from these Chevalley generators
as follows
\begin{equation}\label{uq-com1}
\begin{split}
\ele_{c,a}&= \ele_{c,b}\ele_{b,a}-q\ele_{b,a}\ele_{c,b}\,,\\
\ele_{a,c}&= \ele_{a,b}\ele_{b,c}-q^{-1}\ele_{b,c}\ele_{a,b}\,,\qquad a<b<c\,.
\end{split}
\end{equation}

Using generators $\ele_{ab}$ we  define an evaluation homomorphism:
\begin{equation}\label{eval}
\mathcal{E}v_z\sk{\LL^+(u)}= \LL^+-\, \frac{z}{u}\ \LL^-\,,\quad
\mathcal{E}v_z\sk{\LL^-(u)}= \LL^--\, \frac{u}{z}\ \LL^+\,.
\end{equation}
One can check that the relations \r{uq-com} and \r{uq-com1}
follow from \r{L-op-com}.

Let $M_\Lambda$ be a $\Uqqq$-module generated by a vector $v$,
satisfying the conditions $\ele_{a,a}\,v=q^{\weig_a}\,v$ and
$\ele_{a,b}\,v=0$ for $a<b$. Then $v$ is a singular weight vector of
the evaluation $\Uqgln$ module $M_{\lambda}(z)$. Taking into account
a reordering of the factors
\begin{equation}\label{reorder}
\prod_{b=1}^{N-2}\prod_{\ell=s^b_b+1}^{n_b}\lambda_{b+1}(t^b_\ell)=
\prod_{b=2}^{N-2}\prod_{a=1}^{b-1}\prod_{\ell=1}^{s^b_a}\lambda_{a+1}(t^a_{\ell+\ppp^b_a})
\end{equation}
we can present  the off-shell Bethe vector in $M_{\lambda}(z)$ as
\begin{equation}\label{Wt5}
\begin{split}
\ds &{\bf w}^{N-1}_{M_\Lambda(z)}(\bar{t}_{\segg{\bar{n}}})=\ds
\frac{(q-q^{-1})^{\sum_{a=1}^{N-1}n_a}} {\prod_{a=1}^{N-1}
\prod_{\ell=1}^{n_a} t^a_\ell} \sum_{\admis{s}\prec\, \bar{n}}\left(
\prod_{a\leq b}^{N-1}\frac{1}{[s^b_a-s^b_{a-1}]_q!}
\ \ \times\right.\\
&\quad\times\ds \sk{\ \prod_{N-1\geq b\geq 1}^{\longleftarrow}
\sk{\ \prod_{1\leq a\leq b}^{\longrightarrow}
\sk{z\ele_{b+1,a}\ele^{-1}_{b+1,b+1}}^{s^b_a-s^b_{a-1}}}}v\ \times\\
&\quad\times\ds
\tSym_{\ \bar t_{\segg{\bar n}}}\left.\sk{
\prod_{b=2}^{N-1}\prod_{a=1}^{b-1}\prod_{\ell=1}^{s^b_a}
\frac{q^{\weig_{a+1}}t^{a}_{\ell+\ppp^{b}_{a}}   -q^{-\weig_{a+1}}z }
{1-t^{a}_{\ell+\ppp^{b}_{a}}/t^{a+1}_{\ell+\ppp^{b}_{a+1}}}
\prod_{\ell'=1}^{\ell+\ppp^b_{a+1}-1}
\frac{q-q^{-1}t^{a}_{\ell+\ppp^{b}_{a}}/t^{a+1}_{\ell'}}
{1-t^{a}_{\ell+\ppp^{b}_{a}}/t^{a+1}_{\ell'}}
}\right),
\end{split}
\end{equation}
where $\ppp^b_a=s^a_a+\cdots+s^{b-1}_a$.

\section{Relation to Tarasov-Varchenko construction}

The original inductive construction of nested Bethe vectors of
\cite{KR83} was then developed in \cite{VT}, where these vectors
were defined as certain matrix elements of monodromy operators (see
below). In \cite{KPT} a conjecture about the coincidence of the
construction from \cite{VT} and those used in the present paper was
stated. Recently, Tarasov and Varchenko managed to calculate nested
Bethe vectors in evaluation modules \cite{TV3}. In this section we
use the results of \cite{TV3} to prove a variant of the conjecture
of \cite{KPT}.

\def\Rmat{{\mathcal R}}
\def\Kmat{{\bar{\mathcal R}}}

The $R$-matrix, used in \cite{TV3} differs from \rf{UqglN-R}. To
achieve a compatibility of the results, we slightly modify our
construction. Let $\Rmat\in U_q({\mathfrak{b}}^+)\ot
U_q({\mathfrak{b}}^-)$ be the universal $R$-matrix of $\Uqgln$
(with droped factor $q^{-(c\ot d+ d\ot c)/2}$), and
 $\Kmat$  the universal $R$-matrix of the Hopf subalgebra
$U_q(\mathfrak{gl}_N)$, described in the previous section. Define
$\varepsilon_i\in U_q(\mathfrak{gl}_N)$ by the relation
$q^{\varepsilon_i}=\ele_{ii}$. Define the {\it reduced} $R$-matrix
$\Kmat_0$ of $U_q(\mathfrak{gl}_N)$ by the relation
$$\Kmat_0=
\Kmat q^{\sum_{i=1}^N \varepsilon_i\otimes\varepsilon_i }\,.$$
Then $\Kmat_0^{21}$ is a two-cocycle with respect to comultiplication $\Delta$.
{}For any $x\in\Uqgln$ we set $\tilde{\Delta}(x)=(\Kmat_0^{21})^{-1}\Delta(x)\Kmat_0^{21}$.
The universal $R$-matrix $\tilde{\Rmat}$ for the comultiplication $\tilde{\Delta}$
is  $\tilde{\Rmat}=\Kmat_0^{-1}\Rmat \Kmat_0^{21}$. It is an element
of $U_q(\tilde{\mathfrak{b}}^+)\ot U_q(\tilde{\mathfrak{b}}^-)$, where the
new Borel subalgebras
$U_q(\tilde{\mathfrak{b}}^+)$ and $U_q(\tilde{\mathfrak{b}}^-)$ are
 determined by $L$-operators
$\tilde{\LL}^+(z)=(\pi(z)\ot 1)(\tilde{R}^{21})^{-1})$ and
$\tilde{\LL}^-(z)=(\pi(z)\ot 1)(\tilde{R})$. Here $\pi(z)$ is a
vector representation of $\Uqgln$ evaluated at the point $z$. We have by
the construction
\begin{equation}\label{newL}
\tilde{\LL}^\pm(z)=\left(\LL_0^-\right)^{-1}
\LL^\pm(z)\left(\LL_0^+\right)^{-1}\,,
\end{equation}
where $\LL^\pm(z)$ are $L$-operators of Section \ref{section2.2}, and
$\LL^\pm_0$ are given by the relations \rf{LL+} and \rf{LL-}. The
$L$-operators $\tilde{\LL}^\pm(z)$ satisfy the Yang-Baxter relations with
$R$-matrix \rf{UqglN-R2}, initial conditions \rf{Linit2}
and comultiplication rule \rf{coprL2}:
\begin{equation}
\begin{split}\label{UqglN-R2}
\tilde{\RR}(u,v)\ =\ &\ \sum_{1\leq i\leq N}\E_{ii}\ot \E_{ii}\ +\
\frac{u-v}{qu-q^{-1}v}
\sum_{1\leq i<j\leq N}(\E_{ii}\ot \E_{jj}+\E_{jj}\ot \E_{ii})
\\
+\ &\frac{q-q^{-1}}{qu-q^{-1}v}\sum_{1\leq i<j\leq N}
(u \E_{ij}\ot \E_{ji}+ v \E_{ji}\ot \E_{ij})\,,
\end{split}
\end{equation}
\begin{align}
\label{Linit2}
\tilde{\LL}^{+}_{ji}[0]=\tilde{\LL}^{-}_{ij}[0]=0,\qquad
&\tilde{\LL}^{+}_{kk}[0]\tilde{\LL}^{-}_{kk}[0]=1,
\qquad 1\leq i<j \leq N,
\quad 1\leq k\leq N\, ,\\
\label{coprL2}
&\tilde{\Delta} \sk{\tilde{\LL}^{\pm}_{ij}(u)}=\sum\nolimits_{k}\
\tilde{\LL}^{\pm}_{kj}(u)\otimes
\tilde{\LL}^{\pm}_{ik}(u)\,.
\end{align}
 Let
\begin{equation}\label{L-op2}
\tilde{\LL}^\pm(z)=\sk{\sum_{i=1}^N
\E_{ii}+\sum^N_{i<j}\tilde{\FF}^{\pm}_{j,i}(z)\E_{ij}}
\cdot
\sum^N_{i=1}k^\pm_{i}(z)\E_{ii}
\cdot \sk{\sum_{i=1}^N
\E_{ii}+\sum^N_{i<j}\tilde{\EE}^{\pm}_{i,j}(z)\E_{ji}}\,.
\end{equation}
be the Gauss decomposition of $\tilde{\LL}^\pm(z)$.
We have the following connection to the current realization of $\Uqgln$:
\begin{align}\label{DF-iso2}
&E_i(z)=\tilde{\EE}^{+}_{i,i+1}(z)-\tilde{\EE}^{-}_{i,i+1}(z)\,,
&&F_i(z)=\tilde{\FF}^{+}_{i+1,i}(z)-\tilde{\FF}^{-}_{i+1,i}(z)\,,\\
\label{DFinverse2}
&\tilde{\EE}^{\pm}_{i,i+1}(z)=E_i(z)^{(\pm)}\,,
&&\tilde{\FF}^{\pm}_{i+1,i}(z)=z\left(z^{-1}F_i(z)\right)^{(\pm)}\,.
\end{align}
Besides, the correspondence \rf{newL} imply the relations
\begin{equation}\label{newF}
\tilde{\EE}^{\pm}_{i,i+1}(z)={\EE}^{\pm}_{i,i+1}(z)-E_i[0],\qquad
\tilde{\FF}^{\pm}_{i,i+1}(z)={\FF}^{\pm}_{i,i+1}(z)+F_i[0].
\end{equation}
We have the new decomposition $U_F=\tilde{U}_f^-\tilde{U}_F^+$, where
$\tilde{U}_f^-=U'_F\cap U_q(\tilde{\mathfrak{b}}^-)\,,$ and
$\tilde{U}_F^+=U_F\cap U_q(\tilde{\mathfrak{b}}^+)\,$,
and the new projections
$\tilde{\Pfp}:U_F \to \tilde{U}_F^+$ and
$\tilde{\Pfm}:U_F\to \tilde{U}_f^-$ defined as
\begin{equation}\label{pgln2}
\begin{split}
\tilde{\Pfp}(\Ff_-\ \Ff_+)=\coun(\Ff_-)\ \Ff_+, \qquad
\tilde{\Pfm}(\Ff_-\ \Ff_+)=\Ff_-\ \coun(\Ff_+),
\qquad \Ff_-\in \tilde{U}_f^-,
\quad \Ff_+\in \tilde{U}_F^+ .
\end{split}
\end{equation}
For any $\Pi$-ordered multiset $I$ we set
\begin{equation}\label{W-un-any2}
\tilde{\mathcal{W}}_I(t_i|_{i\in I})=\tilde{\Pfp}\sk{F_{\i(i_n)}(t_{i_n})
F_{\i(i_{n-1})}(t_{i_{n-1}})\cdots F_{\i(i_2)}(t_{i_2}) F_{\i(i_1)}(t_{i_1})}\,.
\end{equation}
A small modification of arguments of
 \cite{KPT} shows that a collection of $V$-valued functions
 $\tilde{w}_{V,I}(t_i|_{i\in I})=\tilde{\mathcal{W}}_I(t_i|_{i\in I})\ v\ $,
  where $V$ is a
$\Uqgln$-module with a singular weight vector $v$,
form a $q$-symmetric weight function, that is, satisfy the setting of Section
\ref{section3} with respect to the comultiplication $\tilde{\Delta}$.

Denote by $\tilde{\mathcal{W}}^{N-1}(\bar{t}_{\segg{\bar{n}}})$
the universal weight function associated with the set of  variables
\r{set111}.
\begin{equation}\label{tuwf}
\tilde{\mathcal{W}}^{N-1}(\bar{t}_{\segg{\bar{n}}})=\tilde{P}^+\left(
F_{N-1}(t_{n_{N-1}}^{N-1})\cdots F_{N-1}(t_{1}^{N-1})
  \quad \cdots\quad
F_{1}(t_{n_1}^1)\cdots F_{1}(t_{1}^1)
\right),
\end{equation}
and by $\tilde{{\bf w}}_V^{N-1}(\bar t_{\segg{\bar n}}) = \qsym(\bar
t_{\bar
n})\prod_{a=2}^{N}\prod_{\ell=1}^{n_{a-1}}\lambda_a(t^{a-1}_\ell) \
\tilde{\mathcal{W}}^{N-1}(\bar t_{\segg{\bar n}})v$ the related modified weight
function, associated with a weight singular vector $v$ of
a $\Uqgln$-module $V$.

Set $\tY(\bar u;\bar v)=\prod_{m=1}^k\frac{v_m}{u_m}\,Y(\bar u;\bar v)$,
 where  the series $Y(\bar u;\bar v)$ is defined by  \r{rat-Y} and let
$\tX(\bar t_{\seg{\bar\rr}{\bar\rr-\bar\lll}})$ be the series defined by
\r{rat-X} with $Y(\cdot)$ replaced by $\tY(\cdot)$.
The first step of the calculation of recurrence relation for  \r{tuwf}, as well
as the definition of strings is unchanged. The difference appears during the calculation
 of the projection \r{rr2a}.
This difference results that the series $X(\bar t_{\seg{\bar n}{\bar n-\bar s}})$
in \r{main-rec} is replaced by $\tX(\bar t_{\seg{\bar n}{\bar n-\bar s}})$
and the recurrence relation for \r{tuwf} takes the form
\begin{equation*}
\begin{split}
\tcW^{N-1}(\bar{t}_{\segg{\bar{n}}})&=
\sum_{\bar s}
\frac{1}{s_1!}\prod_{a=2}^{N-1}\frac{1}{(s_a-s_{a-1})!}
\prod_{a=1}^{N-2}\frac{1}{(n_a-s_{a})!}\times\\[5mm]
&\times\ \qSym_{\ \bar t_{\segg{\bar n}}}
\sk{Z_{\bar\ss}(\bar t_{\segg{\bar n}})\cdot \tX(\bar t_{\seg{\bar
n}{\bar n-\bar s}})\cdot \tilde{\Pfp}\sk{\F^{N-1}_{\bar s}(\bar
t^{N-1}_{\segg{n_{N-1}}})} \cdot \tcW^{N-2}(\bar{t}_{\segg{\bar
n-\bar s}})}.
\end{split}
\end{equation*}
The iteration of this relation  yields an analog of  Theorem~\ref{main-th}
\begin{equation}\label{th1a}
 \tcW^{N-1}(\bar{t}_{\segg{\bar{n}}})=  \qSym
_{\ \bar t_{\segg{\bar n}}}
\sum_{\admis{s}\prec\, \bar{n}}\left(
\frac{\tcZ_{\admis{s}}(\bar t_{\segg{\bar n}})
}{\prod_{a\leq b}(s^b_a-s^b_{a-1})!}
\prod^{\longleftarrow}_{N-1\geq j\geq 1}
\tPfp\sk{\F_{\bar s^j}^j (\bar t^j_{\segg{s^j_j}})    }
\right),
\end{equation}
where
\begin{equation*}
\tcZ_{\admis{s}}(\bar t_{\segg{\bar n}})=
\prod_{b=2}^{N-1}\prod_{a=1}^{b-1}\prod_{\ell=1}^{s^b_a}
\frac{t^{a}_{\ell+\ppp^{b}_{a}}/t^{a+1}_{\ell+\ppp^{b}_{a+1}}}
{1-t^{a}_{\ell+\ppp^{b}_{a}}/t^{a+1}_{\ell+\ppp^{b}_{a+1}}}
\prod_{\ell'=1}^{\ell+\ppp^b_{a+1}-1}
\frac{q-q^{-1}t^{a}_{\ell+\ppp^{b}_{a}}/t^{a+1}_{\ell'}}
{1-t^{a}_{\ell+\ppp^{b}_{a}}/t^{a+1}_{\ell'}}\,,
\end{equation*}
and $\ppp^b_a=s^a_a+\cdots+s^{b-1}_a$.
 An analog
 of  Proposition~\ref{idenGC},
$$\tilde{\Pfp}\sk{F_{j,i}(t)}\,=\,(q-q^{-1})^{j-i-1}\tilde{\FF}^{+}_{j,i}(t)\,,\qquad i<j-1\,.
$$
allows to rewrite \rf{th1a} in Gauss coordinates of $L$-operator $\tilde{\LL}^+(t)$:
\begin{equation*}
\begin{split}
\ds &\tbfw_V^{N-1}(\bar{t}_{\segg{\bar{n}}})=\ds
\qsym(\bar t_{\segg{\bar n}})\
\qSym
_{\ \bar t_{\segg{\bar n}}}
\sum_{\admis{s}\prec\, \bar{n}}\left(
\frac{(q-q^{-1})^{\sum_{b=1}^{N-1}(n_b-s_b^b)}}
{ \prod_{a\leq b}^{N-1}(s^b_a-s^b_{a-1})!}
\ \tcZ_{\admis{s}}(\bar t_{\segg{\bar n}}) \times  \right.\\
 &\left.\prod_{N-1\geq b\geq 1}^{\longleftarrow}
\left(
\prod_{1\leq a\leq b}^{\longrightarrow}
\prod_{\ell=s^b_{a-1}+1}^{s^b_a} \left( \tilde{\FF}^{+}_{b+1,a}(t^b_\ell)k^+_{b+1}(t^b_\ell)
\prod_{\ell'=\ell+1}^{s^b_a}
\frac{t^b_\ell-t^b_{\ell'}}{q^{-1}t^b_\ell-qt^b_{\ell'}} \right)\right)
 \prod_{\ell=s^b_b+1}^{n_b}k^+_{b+1}(t^b_\ell)\right)v.
\end{split}
\end{equation*}
The arguments for the derivation of an analog of
 Theorem~\ref{th-Wt3} and Corollary~\ref{cor53} are unchanged. We get finally
the following expression for the modified weight function:
\begin{equation}\label{Wt444}
\begin{split}
\ds \tbfw_V^{N-1}(\bar{t}_{\segg{\bar{n}}})&=\ds  \tSym
_{\ \bar t_{\segg{\bar n}}}
\sum_{\admis{s}\prec\, \bar{n}}\left((q-q^{-1})^{\sum_{b=1}^{N-1}(n_b-s_b^b)}
\prod_{a\leq b}^{N-1}
\frac{1}{[s^b_a-s^b_{a-1}]_q!}\  \ \times \right.\\
&\times\
\prod_{b=2}^{N-1}\prod_{a=1}^{b-1}\prod_{\ell=1}^{s^b_a}
\frac{t^{a}_{\ell+\ppp^{b}_{a}}/t^{a+1}_{\ell+\ppp^{b}_{a+1}}}
{1-t^{a}_{\ell+\ppp^{b}_{a}}/t^{a+1}_{\ell+\ppp^{b}_{a+1}}}
\prod_{\ell'=1}^{\ell+\ppp^b_{a+1}-1}
\frac{q-q^{-1}t^{a}_{\ell+\ppp^{b}_{a}}/t^{a+1}_{\ell'}}
{1-t^{a}_{\ell+\ppp^{b}_{a}}/t^{a+1}_{\ell'}}
\\
&\times\ds \left.
\prod_{N-1\geq b\geq 1}^{\longleftarrow} \sk{\prod_{1\leq a\leq b}^{\longrightarrow}
\sk{\prod_{\ell=s^b_{a-1}+1}^{s^b_a} \tLL^+_{a,b+1}(t^b_\ell) }
\prod_{\ell=s^b_b+1}^{n_b}\tLL^+_{b+1,b+1}(t^b_\ell)}
\right)v\,.
\end{split}
\end{equation}
The embedded algebra $U_q(\mathfrak{gl}_N)$ is given again by zero modes
of $L$-operators. We set $\tLL^\pm\equiv \tLL^\pm[0]$. Due to \rf{newL}, we have
$\tLL^+=(\LL^-_0)^{-1}{\mathrm K}$ and $\tLL^-={\mathrm K}^{-1}(\LL^+_0)^{-1}$.
Introduce generators $\tele_{i,j}$ of $U_q(\mathfrak{gl}_N)$ by the relations
\begin{equation*}
\tLL^+=(\LL^-_0)^{-1}{\mathrm K}=
\sk{ \begin{array}{cccc}
1&\hba \tele_{21}&\cdots&\hba \tele_{N1}\\
&\ddots&\ddots&\vdots \\
&&1& \hba \tele_{N,N-1}\\
&&&1\\
\end{array} }\sk{ \begin{array}{cccc} \tele_{11} &&&\\
&\tele_{22}&0&\\ & 0&\ddots& \\ &&&    \tele_{NN}
\end{array} }
\end{equation*}
\begin{equation*}
\tLL^-={\mathrm K}^{-1}(\LL^+_0)^{-1}=
\sk{ \begin{array}{cccc} \tele^{-1}_{11} &&&\\
&\tele^{-1}_{22}&0&\\ & 0&\ddots& \\   &&&    \tele^{-1}_{NN}
\end{array} }
\sk{ \begin{array}{cccc}
1&&&\\ -\hba \tele_{12}&1&0&\\ \vdots & \ddots &\ddots&\\
-\hba \tele_{1N}&\cdots&-\hba \tele_{N-1,N}&1
\end{array} }
\end{equation*}
In particular, $\tele_{ij}=\ele_{ij}$ if $|i-j|\leq 1$, so
the Chevalley generators
$\tele_{a,a}$, $\tele_{a,a+1}$ and $\tele_{a+1,a}$ satisfy the same relations
\r{uq-com} and \rf{serregl}
with  different rules for the composed roots generators
\begin{equation*}
\begin{split}
\tele_{c,a}&= \tele_{c,b}\tele_{b,a}-q^{-1}\tele_{b,a}\tele_{c,b}\,,\\
\tele_{a,c}&= \tele_{a,b}\tele_{b,c}-q\tele_{b,c}\tele_{a,b}\,,\qquad a<b<c\,.
\end{split}
\end{equation*}
Define the evaluation homomorphism of  the algebra $\Uqgln$ to $\Uqqq$
as in \r{eval}:
\begin{equation}\label{eval2}
\mathcal{E}v_z\sk{\tilde{\LL}^+(u)}= \tilde{\LL}^+-\, \frac{z}{u}\ \tilde{\LL}^-\,,\quad
\mathcal{E}v_z\sk{\tilde{\LL}^-(u)}= \tilde{\LL}^--\, \frac{u}{z}\ \tilde{\LL}^+\,.
\end{equation}
 The substitution of \rf{eval2} into  \rf{Wt444} gives
\begin{equation}\label{Wt555}
\begin{split}
&\ds \tbfw^{N-1}_{M(z)}(\bar{t}_{\segg{\bar{n}}})=
(q-q^{-1})^{\sum_{a=1}^{N-1}n_a}
\sum_{\admis{s}\prec\, \bar{n}}\left(\sk{
\prod_{N-1\geq b>a\geq 1}^{\longleftarrow}
\frac{q^{s^b_{a-1}(s^b_{a-1}-s^b_{a})}}{[s^b_a-s^b_{a-1}]_q!}
\ \check\ele_{b+1,a}^{s^b_a-s^b_{a-1}}}v\ \times\right.\\
&\quad\times\ds\left.
\tSym_{\ \bar t_{\segg{\bar n}}}\sk{
\prod_{b=2}^{N-1}\prod_{a=1}^{b-1}\prod_{\ell=1}^{s^b_a}
\frac{q^{\weig_{a+1}}t^{a}_{\ell+\ppp^{b}_{a}}   -q^{-\weig_{a+1}}z }
{t^{a+1}_{\ell+\ppp^{b}_{a+1}}-t^{a}_{\ell+\ppp^{b}_{a}}}
\prod_{\ell'=1}^{\ell+\ppp^b_{a+1}-1}
\frac{qt^{a+1}_{\ell'}-q^{-1}t^{a}_{\ell+\ppp^{b}_{a}}}
{t^{a+1}_{\ell'}-t^{a}_{\ell+\ppp^{b}_{a}}}
}\right),
\end{split}
\end{equation}
where $\check\ele_{b+1,a}=\tele_{b+1,a}\tele_{b+1,b+1}$. The factor
$q^{s^b_{a-1}(s^b_{a-1}-s^b_{a})}$ appears after the reordering of the generators
$\check\ele_{b+1,a}$. To obtain \r{Wt555} from \r{Wt444} we use again the
reordering \r{reorder}.

\bigskip

 Remind the construction of
 \cite{VT}. Let $L$-operator $\LL^{}(z)=\sum_{k=0}^\infty
\sum_{i,j=1}^N \E_{ij}\otimes \LL^{}_{ij}[ k]z^{-k}$ of some Borel
subalgebra ${U}_q(\check{\mathfrak{b}}^+)$ of $\Uqgln$ satisfies the Yang-Baxter relation with a
$R$-matrix $\check{\RR}(u,v)$.  We use the notation $\LL^{(k)}(z)\in \sk{\CC^N}^{\otimes
M}\ot {U}_q(\check{\mathfrak{b}}^+)$ for an
$\LL$-operator acting nontrivially on $k$-th
tensor factor in the product $\sk{\CC^{N}}^{\otimes M}$ for $1\leq k\leq M$.
 Consider a series on $M$ variables
\begin{equation}\label{product}
\ct(u_1,\ldots,u_M)=\LL^{(1)}(u_1)\cdots \LL^{(M)}(u_M)
\cdot \ccR^{(M,\ldots,1)}(u_M,\ldots,u_1)
\end{equation}
with coefficients in $\sk{\End(\CC^N)}^{\ot M}\ot {U}_q(\check{\mathfrak{b}}^+)$,
where
\begin{equation}\label{Rproduct}
\ccR^{(M,\ldots,1)}(u_M,\ldots,u_1)=
\prodl_{1\leq i<j\leq M}\check{\RR}^{(ji)}(u_j,u_i)\,.
\end{equation}
In the ordered product of $R$-matrices \r{Rproduct} the factor $\check{\RR}^{(ji)}$ is to
the left of the factor $\check{\RR}^{(ml)}$ if $j>m$, or $j=m$ and $i>l$.
Consider the set of  variables \r{set111}.
Following \cite{VT}, set
\begin{equation}\label{Vit-el}
\bbb(\bar t_{\segg{\bar n}})
={\rm tr}_{\CC^{\ot|{\bar n}|}}\ot {\rm id}\ (\ct(t^1_1,\ldots,t^1_{n_1};\ \ldots\ ;
t^{N-1}_1\!\!,\ldots,t^{N-1}_{n_{N-1}})\E_{21}^{\ot n_1}\ot\cdots\ot \E_{N,N-1}^{\ot
n_{N-1}}\ot 1),
\end{equation}
where $|\bar n|=n_{N-1}+\cdots+n_1$. The element \r{Vit-el} is given by \r{product}
with the identification: $M=|\bar n|$, for $a=1,\ldots,N-1$ and $n_1+\cdots+n_{a-1}<i\leq
n_1+\cdots+n_{a}$, $u_i=t^a_{i-n_1-\cdots-n_{a-1}}$.
The coefficients of
$\bbb(\bar t_{\bar n})$ are elements of the Borel subalgebra
$U_q(\check{\mathfrak{b}}^+)$. For any $\Uqgln$-module $V$ with a
 singular vector $v$ denote
\begin{equation}\label{Vel-vs-pr}
{\mathbb{B}}_V(\bar t_{\segg{\bar n}})=\bbb(\bar t_{\segg{\bar n}})v.
\end{equation}
We can establish a correspondence of the calculations above and of \cite{TV3}
as follows. We identify the $R$-matrix $\tilde{\RR}(u,v)$, see \rf{UqglN-R2} and
$\check{\RR}(u,v)$, the Borel subalgebra $U_q(\tilde{\mathfrak{b}}^+)$ (see the
beginning of this section) with $U_q(\check{\mathfrak{b}}^+)$ of \cite{TV3}.
Set
\begin{equation*}
\eta(\bar t_{\segg{\bar n}})=\prod_{1\leq a<b\leq N-1}\ \
\prod_{1\leq j\leq n_b}\  \ \prod_{1\leq i\leq n_a}\
\frac{qt^b_j-q^{-1}t^a_i}{t^b_j-t^a_i}
\end{equation*}

 Then we have
\begin{equation}\label{tozhd}
\tilde{\bf w}^{N-1}_{M_\Lambda(z)}(\bar t_{\segg{\bar n}})=\eta(\bar t_{\segg{\bar n}})
{\mathbb{B}}_{M_\Lambda(z)}(\bar t_{\segg{\bar n}})
\end{equation}
for arbitrary evaluation module $M_\Lambda(z)$. This is just a literal
coincidence of eq. \rf{Wt555} and  of eq. (6.11) in \cite{TV3}.
Next, we know that the comultiplication property of the weight functions
$\tilde{\bf w}^{N-1}_{V}(\bar t_{\segg{\bar n}})$ and
${\mathbb{B}}_{V}(\bar t_{\segg{\bar n}})$ coincide \cite{KPT}.
Thus
\begin{equation}\label{tozhd2}
\tilde{\bf w}^{N-1}_{1_{g(z_0)}\ot M_{\Lambda_1}(z_1)\ot ...\ot
M_{\Lambda_n}(z_n)}(\bar t_{\segg{\bar n}})=\eta(\bar t_{\segg{\bar n}})
{\mathbb{B}}_{1_{g(z_0)}\ot M_{\Lambda_1}(z_1)\ot ...\ot
M_{\Lambda_n}(z_n)}(\bar t_{\segg{\bar n}})
\end{equation}
for any tensor product of evaluation modules and a one-dimensional
module $1_{g(z_0)}$, in which every $L_{ii}(z)$ acts by
multiplication on $g(z_0)$ (both weight functions in consideration
are trivial for one-dimensional modules). The identity \rf{tozhd2}
is sufficient to conclude the general coincidence of the two
constructions.

Let $J$ be the left ideal of $U_q(\tilde{\mathfrak{b}}^+)$, generated by all
$E_i[n]$, $i=1,\ldots,N-1$, $n>0$ (equivalently, by modes of $\tilde{\E}_{b,a}(t)$,
$a>b$). Let  $\tilde{\mathcal{W}}^{N-1}(\bar{t}_{\bar{n}})$ be the
universal weight function given by eq. \rf{tuwf}.

The following result verifies the conjecture of \cite{KPT} for
the universal weight function \rf{W-un-any2}, and the weight function of
\rf{Vit-el}, related to $R$-matrix \rf{UqglN-R2}.
\begin{theorem}${}$

\noindent
(i) The two  weight functions are equal for each irreducible
finite-di\-men\-si\-o\-nal $\Uqgln$-module $V$ with a singular vector $v$:
\begin{equation*}
\tilde{\bf w}^{N-1}_{V}(\bar t_{\segg{\bar n}})=\eta(\bar t_{\segg{\bar n}})
{\mathbb{B}}_{V}(\bar t_{\segg{\bar n}})\,.
\end{equation*}

\noindent (ii) Consider $U_q(\tilde{\mathfrak{b}}^+)$ as an algebra over
$\CC[[(q-1)]]$. Then
\begin{equation} \label{coin}
\tilde{\mathcal{W}}^{N-1}(\bar{t}_{\segg{\bar{n}}})
\prod_{a=1}^{N-1}
\prod_{\ell=1}^{n_a}\ k^+_{a+1}(t^{a}_\ell)\equiv\eta(\bar t_{\segg{\bar n}})
\bbb(\bar t_{\segg{\bar n}}) \quad\mod J\,.
\end{equation}
\end{theorem}

\noindent {\it Proof}. For the proof of (i) we apply
\rf{tozhd2} and the classical result \cite{CP}: every irreducible
finite-dimensional $\Uqgln$-module with a singular vector $v$ is
isomorphic to a subquotient of a tensor product of one-dimensional
modules and of evaluation modules. More precisely, this subquotient is a
quotient of the submodule of that tensor product  generated by the tensor product
of  weight singular vectors. The singular vector corresponds to
the image of the tensor product of the singular vectors within this
isomorphism.

Due to \rf{tozhd2}, for the proof of (ii) it is sufficient to
verify the following statement. Let $x\in U(\tilde{b}_+)$ be an element of the universal
enveloping algebra of a Borel subalgebra of $\widehat{\mathfrak{gl}}_N$,
which does not belong to $J$. Then there exists a
tensor product of evaluation modules with a weight singular vector $v$ such that
 $xv\not=0$.
This can be observed as follows.  By PBW theorem, we can present
$x$ as a sum of ordered monomials on $\E_{j,i}[n]\in \widehat{\mathfrak{gl}}_N$ with $1\leq i\leq j\leq N$ and
$n\geq 0$. Clearly, $x$ admits a weight decomposition with respect to Cartan
subalgebra of $\mathfrak{gl}_N$. Take a maximal weight component $x_0$.
Present $x_0$ in a form
\begin{equation}\label{pol}
x_0=\sum_k X_kP_k(\{\E_{a,a}[n_a]\}),\
\end{equation}
where each $X_k$ is a monimial as
$X_k=\E_{2,1}[n_1]\E_{2,1}[n_2]\cdots\E_{3,1}[n_l]\cdots \E_{N,N-1}[n_m]$
of a given weight and $P_k(\{\E_{a,a}[n_a]\})$ is a polynomial over commuting
imaginary root vectors $\E_{a,a}[n_a]$, $a=1,\ldots,N$.  Let $m_{i,j}$ denotes the
number of occurrences of all possible $\E_{j,i}[r]$ in the decomposition of $X_1$.

Let $M_{i}$ be a $\mathfrak{gl}_N$ Verma module with a highest vector $v$, satisfying
the relations $\E_{j,i}[0]v=0$ for $i>j$ and $\E_{a,a}[0]v=\delta_{a,i}v$. Let
$M_i(z)$ be the corresponding evaluation module.
For a collection $\bar n=\{n_1,\ldots,n_N\}$ and a set of variables
$z_{\bar n}=\{z_1^1,\ldots,z_{n_1}^1,\ldots,$ $z_{n_N}^N\}$ define $M(z_{\bar n})$
as a tensor product $M(z_{\bar n})=M_1(z_1^1)\ot\cdots\ot M_N(z_{n_N}^N)$.
Its weight singular vector is  $\bar v=v\ot\cdots\ot v$. Denote by
$v_{\{m_{ij}\}}\in M(z_{\bar n})$ the vector
${\E_{2,1}v\ot\cdots\ot \E_{2,1}v}\ot\cdots \ot\E_{N,N-1}v\ot v\ot\cdots\ot v$
where the number of occurrences of each $\E_{i,j}v$ is $m_{i,j}$ (we suppose that
the numbers $n_k$ are big enough). Equip the dual space $M^*(z_{\bar n})$ with a structure
of contragredient $\mathfrak{gl}_N$ module ($\E_{j,i}^*=\E_{i,j}$). Then the vector
$v^*_{\{m_{ij}\}}\in M^*(z_{\bar n})$ is well defined.

Consider the matrix element
$\langle v^*_{\{m_{ij}\}},xv\rangle $.
 By weight arguments it is equal to
$\langle v^*_{\{m_{ij}\}},x_0v\rangle$. It is also clear that we have nonzero
contribution only for the terms with the same number of occurrences of $\E_{j,i}[r]$
(neglecting degrees $r$) as in $X_1$.
In this matrix coefficient the imaginary root vector $\E_{a,a}[m]$ contributes
as the Newton  polynomial $s^a_m=(z^a_1)^m+\cdots+(z^a_{n_a})^m$ of degree $m$ over the variables
$z^a_1,\ldots,z^a_{n_a}$ and thus each polynomial $P_k(\{\E_{a,a}[m_a]\})$ contributes
as this very polynomial over symmetric functions $s^a_m$.

In order to describe the contribution of monomials on $\E_{j,i}$ with
$i<j$ we make the following renaming of the variables $z_i^a$. We
write them first in a natural order $z_1^1,\ldots,z_{n_N}^N$ and then
rename first $m_{1,2}$ of them as $x_{1}^{12},\ldots,x_{m_{12}}^{12}$,
the next $m_{1,3}$ of them as $x_{1}^{13},\ldots,x_{m_{13}}^{13}$ and
so on. The rest of the variables remain unchanged. Then  each
$\E_{j,i}[m]$ with $i<j$ contributes as the Newton polynomial
$s^{ij}_m=(x^{ij}_1)^m+\cdots+(x^{ij}_{m_{ij}})^m$ of degree $m$
 over the variables $x^{ij}_1,\ldots,x^{ij}_{m_{ij}}$.

Now we use the following properties of symmetric (with respect to a
product of symmetric groups) functions: (i) the Newton polynomials of
degree less then the number of the variables are algebraically
independent;
(ii) the products $s_{m_1}s_{m_2}\cdots s_{m_n}$ of the Newton polynomials
$s_k=x_1^k+\ldots +x_n^k$ of $n$ variables taken over all unordered
collections $(m_1,...,m_n)$ of nonnegative integers form a linear basis of the ring of symmetric
functions on $n$ variables;
(iii) for any finite set of linearly independent polynomials of $n$ variables
one can find  $N>n$ such that
all these polynomials  will be linear independent over the ring of symmetric
functions of $N$ variables.

These properties of symmetric functions imply that once we take a
collection $\bar n$ with all $n_i^a$ big enough, the matrix element
$\langle v^*_{\{m_{ij}\}},xv\rangle\, $ is nonzero polynomial over $\{z_{\bar n}\}$.
\hfill{$\square$}
\smallskip

Clearly, the formality restriction in (ii) is artificial and is taken for technical
simplification of the proof.


\setcounter{section}{0}
\setcounter{subsection}{0}

\renewcommand{\thesection}{\Alph{section}}


\section{Analytical properties of composed currents}\label{anal-pr1}

An analytical reformulation of Serre relations \r{serre}, see
\cite{E}, imply the following statements in
the completed algebra $\overline U_F$:

\begin{itemize}
\item[(i)] the products
$\ff_{i}(z)\ff_{i}(w)$  have a simple zero at  $z=w$\footnote{Let $v$ be
a vector of a highest weight module $V$ and $\xi\in V^*$.  The analytical property (i)
means that the  matrix element $\<\xi,\ff_{i}(z)\ff_{i}(w)v\>$
is a meromorphic function of $w$ and $z$ equal to zero at $w=z$.};
\item[(ii)]  the product
$(q^{-1}z_1-q z_2)(z_2-z_3)(qz_1-q^{-1}z_3)\ff_{i-1}(z_1)\ff_{i}(z_2)\ff_{i-1}(z_3)$
vanish on the lines $z_2=z_1=q^{-2}z_3$ and $z_2=q^{-2}z_1=z_3$;
\item[(iii)]  the product
$(z_1-z_2)(q^{-1}z_2-q z_3)(q z_1-q^{-1}z_3)\ff_{i}(z_1)\ff_{i-1}(z_2)\ff_{i}(z_3)$
vanish on the lines $z_2=z_1=q^{2}z_3$ and $z_2=q^{2}z_1=z_3$.
\end{itemize}
The properties (i), (ii) and (iii) imply
 the commutativity
\begin{equation}\label{ser-n2}
[\ff_{i-1}(z)\ff_i(z)\ff_{i+1}(z),\ff_i(w)]=0.
\end{equation}
Let us prove \r{ser-n2}. The basic relations \rf{gln-com} imply that:
\begin{equation}\label{gln-con}
(q^{-1}z-q^{}w)(z-w)(q^{}z-q^{-1}w)[\ff_{i-1}(z)\ff_i(z)\ff_{i+1}(z),\ff_i(w)]=0\,.
\end{equation}
Consider the product $\ff_{i-1}(z)\ff_i(z)\ff_{i+1}(z)\ff_i(w)$.
Due to \r{gln-com} this product as a function
of $w$ has  poles at the points $w=q^{-2}z,z,q^{2}z$. On the other hand it has
zero at $w=z$ due to (i). Due to (ii) the product $\ff_i(z)\ff_{i+1}(z)\ff_i(w)$
has zero at $w=q^{-2}z$ and due to (iii) the product $\ff_{i-1}(z)\ff_i(z)\ff_i(w)$
has zero at $w=q^2z$. This means that the whole product
$\ff_{i-1}(z)\ff_i(z)\ff_{i+1}(z)\ff_i(w)$ has no zeros and no poles. The same is true
for the inverse product $\ff_i(w)\ff_{i-1}(z)\ff_i(z)\ff_{i+1}(z)$.
The relations \r{gln-con} implies now the commutativity
in \r{ser-n2}.

We will use \r{res-in} in order to describe the analytical properties of the product
of composed currents.
\begin{proposition} \label{prop4.2}
The following relations hold in $\overline U_F$ for any $a<b$ and $c<d$:
\begin{eqnarray}
\ff_{b,a}(z)\ff_{d,c}(w)&=&\ff_{d,c}(w)\ff_{b,a}(z)\,,\quad b<c\label{5.62}\\
(q^{-1}z-qw)\ff_{b,a}(z)\ff_{d,c}(w)&=&(z-w)\ff_{d,c}(w)\ff_{b,a}(z)\,,
\quad b=c\label{5.61}\\
\ff_{b,a}(z)\ff_{d,c}(w) &=&  \frac{q-q^{-1}z/w}{1-z/w}\ff_{d,c}(w)\ff_{b,a}(z)\,,
\quad a>c,\ b=d \label{com-aa}\\
\frac{q-q^{-1}w/z}{1-w/z}\ff_{b,a}(z)\ff_{d,c}(w) &=&
 \frac{q-q^{-1}z/w}{1-z/w}\ff_{d,c}(w)\ff_{b,a}(z)\, \quad a=c,\ b=d \label{com-aa1}\\
\ff_{b,a}(z)\ff_{d,c}(w)&=&\ff_{d,c}(w)\ff_{b,a}(z)\,,\quad a<c<d<b \label{cwm1}\\
\frac{q-q^{-1}w/z}{1-w/z}\ff_{b,a}(z)\ff_{d,c}(w)&=& \ff_{d,c}(w)\ff_{b,a}(z)
\,,\quad a=c,\ b<d \label{cwm11}
\end{eqnarray}
\end{proposition}
We note that in the analytic language, both sides of all relations in
\rf{5.62}--\r{cwm11} are analytical functions in $\sk{\CC^*}^2$. This means,
for instance, that the product $\ff_{d,a}(w)\ff_{b,a}(z)$ in \r{cwm11} has no zeroes
and no poles for $b<d$, while the product $\ff_{b,a}(z)\ff_{d,a}(w)$ has a simple
zero at $z=w$ and a simple pole at $z=q^{-2}w$.

Let us prove \r{cwm11}. Consider the product $\ff_{b,a}(z)\ff_{d,a}(w)$.
Due to  \r{ser-n2} and  \r{res-in},
 analytical properties of this product are the same as of
$\ff_a(z)\ff_a(w)\ff_{a+1}(w)$. This product considered as a function of $z$
has poles at the points $z=q^{-2}w$ and $z=q^{2}w$. On the other hand due to
(i) it has a zero at $z=w$ and due to (ii) a zero at $z=q^2w$. It means that the product
$\ff_{b,a}(z)\ff_{d,a}(w)$ has a simple pole at $z=q^{-2}w$ and a simple zero at $z=w$.
Consider now the inverse product $\ff_{d,a}(w)\ff_{b,a}(z)$. Its analytical properties
are defined again due to \r{ser-n2} by the properties of the product $\ff_a(w)\ff_{a+1}(w)\ff_a(z)$.
The latter  has  poles at $z=q^2w$ and at $z=w$. Due to (i) and (ii) it has
zeros at the same points. It means that the inverse product
$\ff_{d,a}(w)\ff_{b,a}(z)$ of composed currents are analytic function in
$\sk{\CC^*}^2$ without zeros and poles. Now the relations \r{gln-com} imply
\r{cwm11}. \hfill$\square$

\begin{corollary}\label{prop46}
The inverse string \r{string-inv} have simple zeros at $t^j_{k}=t^j_m$ and simple poles
at $t^j_{k}=q^{-2}t^j_m$ for all pairs $(m,k)$ such that $\rr_j\geq m> k>\lll_j$.
These poles are the only singularities of the inverse string.
\end{corollary}
{\it Proof}\ \ follows from the analysis of the ordering of the product of currents
in the inverse string. It contains the products of the composed currents
$\ff_{j+1,a}(t^j_k)\ff_{j+1,c}(t^j_m)$, where always $a\leq c$ and $k<m$. Properties
\r{com-aa} and \r{com-aa1} of the products of the composed currents shows that
$\ff_{j+1,a}(t^j_k)\ff_{j+1,c}(t^j_m)$ will always have the pole at $t^j_{k}=q^{-2}t^j_m$
and zero at $t^j_{k}=t^j_m$. \hfill$\square$
\smallskip

The assertion~\ref{prop46} will be used in the next Appendix in order to calculate the projections
of the inverse string \r{string-inv}.

\section{Projections of the string}\label{pr-st3}
In this Appendix we prove Proposition \ref{prop1} describing the
projection $\Pfp$ of the inverse string. We also describe  projection $\Pfm$.
Up to renaming the variables the inverse string is a product
$\ff_{j,i_1}(t_1)\cdots \ff_{j,i_n}(t_n)$ for
$1\leq i_1\leq \cdots \leq i_n  <j$.
\begin{lemma}\label{prop4.6}
For any $1\leq i_1\leq \cdots \leq i_n\leq k  <j$ we have an equality
\begin{equation}\label{RF11}
\begin{split}
\Pfp\!\sk{\ff_{j,i_1}(t_1)\cdots \ff_{j,i_n}(t_n)\ff_{j,k}(t)}&=
\Pfp\!\sk{\ff_{j,i_1}(t_1)\cdots \ff_{j,i_n}(t_n)}\cdot \Pfp\!\sk{\ff_{j,k}(t)}+\\
&+\sum_{b=1}^{n}\frac{V_b(t_1,\ldots,t_{n})} {t-q^2t_b}.
\end{split}
\end{equation}
where $V_b(t_1,\ldots,t_{n})$ take value in $U_f^+$ and do not depend on
the variable $t$.
\end{lemma}
{\it Proof.}\
We want to move the difference $\ff_{j,k}(t)-\Pfp\sk{\ff_{j,k}(t)}$  to the left of the product
$\ff_{j,i_1}(t_1)\cdots \ff_{j,i_n}(t_n)$  and observe that
during this process  appear only simple poles at the points $t=q^2t_b$  with
(independent on $t$)
operator-valued residues.

A particular case of the formula \r{rec-f112a} for $a=i+1$ can be written in the form
\begin{equation}\label{comf-m}
F_{j,k}(t)\,=\,
S_k(F_{j,k+1}(t))+
(q^{-1}-q)F_k(t)^{(-)}\,F_{j,k+1}(t)\,.
\end{equation}
Iterating \r{comf-m} we obtain
\begin{equation}\label{comf-mmm}
F_{j,k}(t)\,=\,
S_k\,S_{k+1}\cdots S_{j-2}(F_{j-1}(t))+\hat F_{j,k}(t)\,,
\end{equation}
where we set  $\hat\ff_{k+1,k}(t)=0$ and
\begin{equation}\label{tFp}
\hat F_{j,k}(t)=(q^{-1}-q)\sum_{m=k+1}^{j-1}
 \sk{\ff_{m,k}(t)^{(-)}-\hat \ff_{m,k}(t)^{(-)}} \ff_{j,m}(t)\,,
\end{equation}
if $j-k>1$.
Proposition \ref{com-cur2} and \r{comf-mmm} imply that
\begin{equation}\label{cur-pre}
\ff_{j,k}(t)-\Pfp\sk{\ff_{j,k}(t)}=\hat F_{j,k}(t)^{(+)}-\ff_{j,k}(t)^{(-)}\,.
\end{equation}
and $\Pfp\sk{\hat F_{j,k}(t)}=0$. Applying Proposition \ref{com-cur2} and
\r{comf-mmm} once again, we get also an equality
 $\Pfp\sk{F_{j,k}(t)^{(-)}}=0$.

Permuting the difference $\ff_{j,k}(t)-\Pfp\sk{\ff_{j,k}(t)}$ and the product
of the currents $\ff_{j,i_1}(t_1)\cdots \ff_{j,i_n}(t_n)$ we  consider the terms
$\hat F_{j,k}(t)^{(+)}$ and $\ff_{j,k}(t)^{(-)}$ in the
r.h.s. of \r{cur-pre} separately.

\begin{lemma}\label{co-min}
For $1\leq i_1\leq\cdots\leq i_n\leq k<j$
\begin{equation}\label{tFnew}
\Pfp\sk{\ff_{j,i_1}(t_1)\cdots \ff_{j,i_n}(t_n)\cdot \hat F_{j,k}(t)}=0\,.
\end{equation}
\end{lemma}
Note that the equality \rf{tFnew} is equivalent to a system of two equalities, each
of them is used further:
\begin{equation}\label{tFnewpm}
\Pfp\sk{\ff_{j,i_1}(t_1)\cdots \ff_{j,i_n}(t_n)\cdot \hat F_{j,k}(t)^{(\pm)}}=0\,.
\end{equation}
{\it Proof}.\ \ We use presentaion \r{tFp}, the equality
$\Pfp\sk{F_{j,k}(t)^{(-)}}=0$ and prove the statement of Lemma by induction
over $k$. For $k=j-2$ the statement \rf{tFnew} follows from the   relations
\begin{equation}\label{cwm}
\begin{split}
\ff_{j,s}(t')\ff_{m,k}(t)^{(-)}&=\ff_{m,k}(t)^{(-)}\ff_{j,s}(t')\,,\quad s<k\\
\ff_{j,k}(t')\ff_{m,k}(t)^{(-)}&=\frac{q^{-1}t'-qt}{t'-t} \ff_{m,k}(t)^{(-)}\ff_{j,k}(t')-
\frac{(q^{-1}-q)t'}{t'-t}\ff_{m,k}(t')^{(-)}\ff_{j,k}(t'),
\end{split}
\end{equation}
valid for $m<j$. The latter relations are the result of applying
 the integral transformation $-\oint \frac{dz}{z}\frac{1}{1-t/z}$ to
 \r{cwm1}, \r{cwm11} with renaming $w\to t'$.
For $k=j-3$ we use again  \r{cwm} and  relations \rf{tFnewpm} for $k=j-2$,
which are already proved and so on.
\hfill$\square$

\begin{lemma}\label{co-min1}
The following relations hold for $i<k$
\begin{equation*}
\begin{split}
\ff_{j,i}(t')\ff_{j,k}(t)^{(-)}&=\frac{t-t'}{q^{-1}t-qt'} \ff_{j,k}(t)^{(-)}\ff_{j,i}(t')+
\frac{(q^{-1}-q)t'}{q^{-1}t-qt'}\ff_{j,i}(t')\ff_{j,k}(t')^{(-)}\,,\\
\ff_{j,k}(t')\ff_{j,k}(t)^{(-)}&=\frac{qt-q^{-1}t'}{q^{-1}t-qt'} \ff_{j,k}(t)^{(-)}\ff_{j,k}(t')+\\
&+
\frac{(q^{-1}-q)t'}{q^{-1}t-qt'}\sk{\ff_{j,k}(t')\ff_{j,k}(t')^{(-)}+\ff_{j,k}(t')^{(-)}
\ff_{j,k}(t')}\,.\\
\end{split}
\end{equation*}
\end{lemma}
\noindent
{\it Proof}.\ Apply
the integral transform $-\oint \frac{dz}{z}\frac{1}{1-t/z}$ to
the relations \r{com-aa}, \r{com-aa1} and rename $w\to t'$.\hfill$\square$

Lemmas~\ref{co-min1} and \ref{co-min} imply Lemma~\ref{prop4.6}.
Exact form of  $V_b(t_1,\ldots,t_n)$ is not important.
\hfill$\square$
\smallskip

Let $\ff_{k,j}(t;t_1,\ldots,t_n)$ be the linear combination of currents introduced by
the relation \r{long-cur}. Lemma~\ref{prop4.6} implies now the following
\begin{lemma}\label{prop710}
For $1\leq i_1\leq \cdots\leq i_n\leq k <j$ we have the equality of formal  series
\begin{equation*}
\Pfp\!\sk{\ff_{j,i_1}(t_1)\cdots \ff_{j,i_n}(t_n)\ff_{j,k}(t)}=
\Pfp\!\sk{\ff_{j,i_1}(t_1)\cdots \ff_{j,i_n}(t_n)}
\Pfp\sk{\ff_{j,k}(t;t_1,\ldots,t_n)}\,.
\end{equation*}
\end{lemma}

{\it Proof}.\
Corollary~\ref{prop46} to Proposition \ref{prop4.2} states that the product
of the currents $\ff_{j,i_1}(t_1)\cdots \ff_{j,i_n}(t_n)\ff_{j,k}(t)$
 has simple zeroes at hyperplanes
$t=t_i$, $i=1,\ldots,n$. Substituting these conditions in \rf{RF11}
 we get a systems of $n$ linear equations over the field of rational
functions $\CC(t_1,\ldots,t_{n})$ for the operators $V_b(t_1,\ldots,t_{n})$:
\begin{equation}\label{system}
\sum_{b=1}^{n}\frac{V_b(t_1,\ldots,t_{n})} {t_a-q^2t_b}=
V\cdot \Pfp\sk{\ff_{k,j}(t_a)}\,, \quad i=1,\ldots, n\, ,
\end{equation}
where $V=\Pfp\sk{\ff_{j,i_1}(t_1)\cdots \ff_{j,i_n}(t_n)}$.
The determinant of the matrix $B_{a,b}=(t_a-q^2t_b)^{-1}$ of this system is nonzero in
 $\CC(t_1,\ldots,t_{n})$,
$$\det(B)=(-q^2)^{\frac{n(n+1)}{2}}\frac{\prod_{a\not=b}(t_a-t_b)^2}
{\prod_{a,b}(t_a-q^2t_b)}\,,$$
 hence the system  has a unique solution over  $\CC(t_1,...,t_{n})$.
This implies that the operators $V_b$ are linear combinations over
$\CC(t_1,...,t_{a-1})$ of the operators $V\cdot \Pfp\sk{\ff_{j,k}(t_a)}$, $a=1,...,n$, and
the projection
\begin{equation}\label{4.21}
\begin{split}
\Pfp\sk{\ff_{j,i_1}(t_1)\cdots \ff_{j,i_n}(t_n)\ff_{j,k}(t)}&=
V\cdot \Pfp\sk{\ff_{j,k}(t)} -\\
&-\sum_{b=1}^{n}\varphi_{t_b}(t;t_1,\ldots,t_{n})\ V\cdot
\Pfp\sk{\ff_{j,k}(t_b)}\,,
\end{split}
\end{equation}
where
$\varphi_{t_b}(t;t_1,\ldots,t_{n})$ $=$
${A_b(t;t_1,\ldots,t_{n})}/
{\prod_{m=1}^{n}(t-q^2t_m)}$ are rational functions  whose
numerators are polynomials on $t$ of degree less then $n$.
The system \rf{system} is satisfied if the rational functions
$\varphi_{t_b}(t;t_1,\ldots,t_{n})$ have the property
$$\varphi_{t_b}(t_a;t_1,\ldots,t_{n})=\delta_{a,b}\,,\qquad a,b=1,...,n\, .$$
This interpolation problem has a unique solution given by formula \r{rat-f}.
This proves Lemma \ref{prop710}. \hfill{$\square$}

The iteration of Lemma \ref{prop710} proves Proposition \ref{prop1}. \hfill{$\square$}

\medskip

Now we describe the projection $\Pfm$ of the composed currents and strings.
Let $\tilde {\rm ad}_x$ denotes the adjoint action in $\Uqgln$ related to coproduct
$\Delta^{op}$,
$
\tilde {\rm ad}_x\ y= \sum\nolimits_l x''_l\cdot y\cdot a^{-1}(x'_l)\ ,$
where
$\Delta x =\sum\nolimits_l x'_l\ot x''_l$.
Define screening operators $\tilde S_i$, $i=1,...,N-1$ by the relation
\begin{equation*}
\tilde S_i\, (y)=\tilde {\rm ad}_{F_{i}[0]}\ (y)=
F_{i}[0]\ y - k^{-1}_i k_{i+1}\ y\  k^{-1}_{i+1} k_i\  F_{i}[0]
\end{equation*}
Inductive using of \r{rec-f112b} implies an analog of Proposition \ref{com-cur2}
 (we do not use it further):
\begin{equation*}
\Pfm\sk{\ff_{j+1,i}(t)}=\tilde S_{j}
                       \tilde S_{{j-1}}
                       \cdots
                       \tilde S_{{i+1}} \sk{ \Pfm( {F}_{i}(t)) }=
                       -\left(\tilde S_{j}
                       \tilde S_{{j-1}}
                       \cdots
                       \tilde S_{{i+1}} \sk{  {F}_{i}(t) }\right)^{(-)}.
\end{equation*}
Set
\begin{equation}\label{rat-fm}
\tilde \varphi_{u_m}(u;u_1,\ldots,u_n)=\prod_{a=1,\ a\neq m}^{n}
\frac{u-u_a}{u_m-u_a}\prod_{a=1}^{n}\frac{qu_m-q^{-1}u_a}{qu-q^{-1}u_a},\qquad
\text{and}
\end{equation}
\begin{equation*}
\tilde{\ff}_{j+1,i}(u;u_{1},\ldots,u_n)={\ff}_{j+1,i}(u)-
\sum_{m=1}^n \tilde{\varphi}_{u_m}(u;u_{1},\ldots,u_n)\ff_{j+1,i}(u_m), \quad 1\leq i\leq j<N.
\end{equation*}
The series $\tilde{\ff}_{i,j+1}(u;u_{1},\ldots,u_n)$ admits another presentation, which
will be used below. Namely, in the notation $u\equiv u_0$ and
 $$\psi_{u_m}(u_0;u_1,\ldots,u_n)=\prod_{j=1}^{n}\dfrac{u_{0}-u_{j}}{qu_{0}-q^{-1}u_j}
\cdot\frac{ \prod_{j=1}^{n}(qu_k-q^{-1}u_j)}
{ \prod_{j=0\ j\neq k}^{n}(u_k-u_j)}$$
we have
\begin{align} \notag
\tilde{\ff}_{j+1,i}(u_{0};u_{1},\ldots,u_n)&=
\sum
_{k=0}^n
\psi_{u_m}(u_0;u_1,\ldots,u_n)
\ff_{j+1,i}(u_k),\qquad \text{so that}\\ \label{pr-st-mm}
\Pfm\sk{\tilde{\ff}_{j+1,i}(u_{0};u_{1},\ldots,u_n)}&=
\sum
_{k=0}^n
\psi_{u_m}(u_0;u_1,\ldots,u_n)
\Pfm\sk{\ff_{j+1,i}(u_k)}.
\end{align}
We have an analog of Proposition \ref{prop1}:
\begin{equation}\label{mpfac}
\begin{split}
\Pfm\!\sk{\ff_{j,i_1}(t_1)\cdots \ff_{j,i_n}(t_n)}&=
\Pfm\!\sk{\tilde{\ff}_{j,i_1}(t_1;t_2,\ldots,t_n)}
\Pfm\!\sk{\tilde{\ff}_{j,i_2}(t_2;t_3,\ldots,t_n)}\cdot\\
&\cdots \Pfm\!\sk{\tilde{\ff}_{j,i_{n-1}}(t_{n-1};t_n)}
\Pfm\!\sk{\ff_{j,i_n}(t_n)}\,.
\end{split}
\end{equation}
Using \r{mpfac} we  present the projection of the string
$\Pfm\sk{\F^{j}_{\bar s'}({\bar t}^{j}_{\segg{s_{j}}})}$ in a factorized form
\begin{equation}\label{pr-st-m}
\begin{split}
\ds \Pfm\sk{\F^j_{\bar s'}(\bar t^j_{\segg{s_{j}}})}&=
\ds \prod^{\longrightarrow}_{1\leq a\leq j}\sk{
\prod^{\longrightarrow}_{\ss_{a-1}< \ell\leq \ss_{a}}
\Pfm\sk{\tilde{F}_{j+1,a}(t^j_\ell;t^j_{\ell+1},\ldots,t^j_{s_j} )}}
\ \times
\\ &\ds\times
\prod_{1\leq\ell<\ell'\leq s_j}
\frac{q^{-1}-qt_\ell^j/t^j_{\ell'}}
{1-t_\ell^j/t^j_{\ell'}}
\prod_{1\leq a\leq j}\sk{
\prod_{s_{a-1}<\ell<\ell'\leq s_a}
\frac{1-t_\ell^j/t^j_{\ell'}}
{q^{}-q^{-1}t_\ell^j/t^j_{\ell'}}}.
\end{split}
\end{equation}

\section{String product expansion}\label{pr6.7}

Let $U_{i,j}$ be subalgebra of $\overline U_F$ generated by the modes of the currents
$\ff_{i}(t),\ldots,\ff_{j}(t)$ and $U^\coun_{i,j}=U_{i,j}\cap{\rm Ker}\coun$ be
 the corresponding augmentation ideal.
\begin{lemma}\label{lemap2} For any $i$ and $j$ with $1\leq i<j\leq N$ we have
\begin{equation}\label{np-vs-m}
\Pfm\sk{F_{j,i}(t)}=-F_{j,i}(t)^{(-)}\quad\mod\quad \Pfm\sk{U^\coun_{i,j-2}}\cdot U_{i,j-1}\,.
\end{equation}
\end{lemma}
{\it Proof}\ \ is based on the relation
\begin{equation*}
\begin{split}
\Pfm\sk{F_{j,i}(t)}+F_{j,i}(t)^{(-)}&=S_{i}\sk{\Pfm\sk{F_{j,i+1}(t)}+F_{j,i+1}(t)^{(-)}}\\
&\mod\ \Pfm\sk{U^\coun_{i,i}}\cdot U_{i,j-1}
\end{split}
\end{equation*}
which is direct consequence of \r{comf-m}. Using this relation several times we obtain
\begin{equation*}
\begin{split}
\Pfm\sk{F_{j,i}(t)}+F_{j,i}(t)^{(-)}&=S_{i}\cdots S_{j-2}\sk{\Pfm\sk{F_{j,j-1}(t)}+F_{j,j-1}(t)^{(-)}}\\
&\mod\ \Pfm\sk{U^\coun_{i,j-2}}\cdot U_{i,j-1}
\end{split}
\end{equation*}
But for the simple roots currents we have
$\Pfm\sk{F_{j,j-1}(t)}+F_{j,j-1}(t)^{(-)}=0$. so Lemma is proved.
Using the commutativity of screening operators and the projections
\cite{KP}, we prove the lemma.
\hfill$\square$
\smallskip

Let $\bar s=\{s_{j+1},s_{j},\ldots,s_2,s_1\}$
 be a collection of non-negative  integers satisfying
admissibility conditions: $s_{j+1}\geq s_{j}\geq \cdots \geq s_1\geq s_0=0$.
Set $\bar s'=\{0,s_{j},\ldots,s_1\}$.
\begin{proposition}\label{main-fact}
For any product  $\F(\bar t^{j+1}_{\segg{s_{j+1}}})$
and  a string $\F^{j}_{\bar s'}({\bar t}^{j}_{\segg{s_{j}}})$
we have an equality
\begin{equation}\label{ca16}
\begin{split}
&\F(\bar t^{j+1}_{\segg{s_{j+1}}}) \cdot\Pfm\sk{
\F^{j}_{\bar s'}({\bar t}^{j}_{\segg{s_{j}}})}=\prod_{i=1}^{j+1}\frac{1}{(s_{i}-s_{i-1})!}\ \
\qSym_{\ t^{j+1}_{\segg{s_{j+1}}}}\left(
\F^{j+1}_{\bar s}({\bar t}^{j+1}_{\segg{s_{j+1}}})\times\right.\\
&\left.\ \times\ \qSym_{\ \bar t^j_{\seg{s_j}{s_{j-1}}}}\cdots
\qSym_{\ \bar t^j_{\seg{s_2}{s_{1}}}}
\qSym_{\ \bar t^j_{\segg{s_1}}}
\sk{ Y(t^{j+1}_{{s_j}},\ldots,  t^{j+1}_{{1}};t^j_{s_j},\ldots,t^j_{1})}\right)
\end{split}
\end{equation}
modulo $\Pfm\sk{U^\coun_{j}}\cdot U_{j+1}$.
\end{proposition}

\noindent
{\it Proof}. \,
Substitute \r{pr-st-m} instead of the second factor of the product
$\F(\bar t^{j+1}_{\segg{s_{j+1}}}) \cdot\Pfm\sk{
\F^{j}_{\bar s'}({\bar t}^{j}_{\segg{s_{j}}})}$. Our strategy is to move each
 factor $\Pfm\!\!\sk{\tilde{F}_{j+1,a}(t^j_\ell;t^j_{\ell+1},\ldots,t^j_{s_j} )}$
to the left of the product
$\F(\bar t^{j+1}_{\segg{s_{j+1}}})$ and keep the terms modulo
$\Pfm(U^\coun_{1,j})\cdot U_{1,j+1}$.
We start from the most left factor taken at $a=1$ and $\ell=1$.
 We  replace first $\Pfm\!\sk{\tilde{F}_{j+1,1}(t^j_1;t^j_{2},\ldots,t^j_{s_j} )}$
 with a linear combination of
$F_{j+1,1}(t^j_{k} )^{(-)}$  with rational coefficients.
This can be done by due to Lemma~\ref{lemap2} and \rf{pr-st-m}, since the modes of the
 $F_{j+1}(t)$ commute with any elements from $\Pfm(U^\coun_{1,j-1})$ and thus
can be moved to the left forming modulo terms. Then we use  the  relations
\rf{com-aa} and ordering rules
\begin{equation}\label{sl3-gen}
  \begin{split}
\ff_{j+1}(t)\ff_{j+1,a}(t')^{(-)}=&-\frac{1}{1-t'/t}\ \ff_{j+2,a}(t)\ +\\
\quad +\  \left(\frac{q-q^{-1}t'/t}{1-t'/t}\right.&\ff_{j+1,a}(t')^{(-)}\left.-\
\frac{q-q^{-1}}{1-t'/t}\ff_{j+1,a}(t)^{(-)}\right)\ff_{j+1}(t),
\end{split}
\end{equation}
for $a=1,\ldots,j$, which are consequence of  \r{5.61}. They give the equalities
\begin{equation}\label{com-1m}
\begin{split}
&F_{j+1}(t^{j+1}_{s_{j+1}})\cdots F_{j+1}(t^{j+1}_{1})\cdot
F_{j+1,1}(t)^{(-)}=
-\sum_{a=1}^{s_{j+1}}\ \prod_{a<\ell\leq s_{j+1}}
\frac{q^{-1}-q^{}t^{j+1}_a/t^{j+1}_\ell}{1-t^{j+1}_a/t^{j+1}_\ell}
\\ &  \prod_{1\leq\ell< a}
\frac{q^{-1}-q^{}t^{j+1}_\ell/t^{j+1}_a}{1-t^{j+1}_\ell/t^{j+1}_a}
 \frac{1}{1-t/t^{j+1}_a}
\ \ F_{j+2,1}(t^{j+1}_a)\
\mathop{\underbrace{F_{j+1}(t^{j+1}_{s_{j+1}})\cdots F_{j+1}(t^{j+1}_{1})}}
\limits_{\
 F_{j+1}(t^{j+1}_a)\ {\rm  omitted}}
\end{split}
\end{equation}
modulo $\Pfm(U^\coun_{1,j})\cdot U_{1,j+1}$. The iteration of \rf{com-1m} using
other ordering rules, being consequences of \r{cwm11} and \rf{cwm1},
\begin{equation*}
\begin{split}
F_{j+2,a}(t)F_{j+1,a}(t')^{(-)}&= \sk{\frac{q^{-1}-qt'/t}{1-t'/t}
F_{j+1,a}(t')^{(-)}- \frac{q^{-1}-q}{1-t'/t}F_{j+1,a}(t)^{(-)}}F_{j+2,a}(t),\\
F_{j+2,a}(t)F_{j+1,b}(t')^{(-)}&=F_{j+1,b}(t')^{(-)}F_{j+2,a}(t).\qquad a<b
\end{split}
\end{equation*}
 leads to the following result.  Denote
\begin{equation}\label{sh-not}
\alpha_1(x)=\frac{q^{-1}-qx}{1-x}\,,\quad
\alpha_2(x)=\frac{q^{-1}-qx}{q-q^{-1}x}\,,\quad
\alpha_3(x)=\frac{q-q^{-1}x}{1-x}\,.
\end{equation}
Then the  product $\F(\bar t^{j+1}_{\segg{s_{j+1}}}) \cdot\Pfm\sk{
\F^{j}_{\bar s'}({\bar t}^{j}_{\segg{s_{j}}})}$ modulo
$\Pfm(U^\coun_{1,j})\cdot U_{1,j+1}$ is equal to the following sum:
\begin{equation}\label{ca11}
\begin{split}
&
\prod_{i=1}^{j}\sk{
\prod_{s_{i-1}<\ell<\ell'\leq s_i}
\alpha_3(t^j_{\ell}/t^j_{\ell'})^{-1} }
\sum_{A_j}\cdots \sum_{A_2}\sum_{A_1}\
Y(t^{j+1}_{a_{s_j}},\ldots,  t^{j+1}_{a_{1}};t^j_{s_j},\ldots,t^j_{1})\\
 &\quad \times   \prod_{i=1}^j
\prod_{\ell<\ell'\atop \ell,\ell'\in A_i}
\!\alpha_1(t^{j+1}_{\ell}/t^{j+1}_{\ell'})
\prod_{\ell<\ell'\atop \ell\in \coA_{i-1}\ \ell'\in A_i}
\!\!\! \alpha_2(t^{j+1}_{\ell}/t^{j+1}_{\ell'})
\prod_{\ell\in \coA_{i-1}\   \ell'\in A_i}
\!\alpha_3(t^{j+1}_{\ell}/t^{j+1}_{\ell'}) \\
&\quad \times
\prod^{\longrightarrow}_{1\leq i\leq j}
\sk{ F_{j+2,i}(t^{j+1}_{a_{s_{i-1}+1}})\cdots F_{j+2,i}(t^{j+1}_{a_{s_i}})}\
\mathop{\underbrace{F_{j+1}(t^{j+1}_{s_{j+1}})\cdots F_{j+1}(t^{j+1}_{1})}}
\limits_{\
 F_{j+1}(t^{j+1}_{a_1}),\ldots,F_{j+1}(t^{j+1}_{a_{s_j}})\ {\rm  omitted}}.
\end{split}
\end{equation}
 The sum in \rf{ca11} goes over all non-ordered subsets $A_1,A_2,\ldots,A_j$ of
the set  $A=\{1,\ldots,s_{j+1} \}$.
These subsets are defined as follows.
Let $A_0=\varnothing$ be an empty set. Denote  $\coA_0=A$. Define inductively
the subsets $A_i$ and $\coA_i$ of $\coA_0$ for
$i=1,\ldots,j$ by the relations: $A_i\cup\coA_i=\coA_{i-1}$. Denote the elements of the set
$A_i$ by the letters $a_{s_{i-1}+1},\ldots,a_{s_i-1},a_{s_i}$. Number of the elements in the
subset $A_i$ is equal to $s_i-s_{i-1}$.
Note also that because
subsets $A_k$ is defined inductively through previous subsets $A_{k-1},\ldots,A_1$ the summations
in \r{ca11} are not commutative. First, we have to sum over all possible subsets $A_1$ in $\coA_0$, then
over all possible $A_2$ in $\coA_1$ and so on.

The last line in \r{ca11} is the inverse string. We can use the  relations \r{com-aa}
between the composed currents and the last product of the rational series in the second line of \r{ca11}
to transform it to the string, again modulo $\Pfm(U^\coun_{1,j})\cdot U_{1,j+1}$:
\begin{equation*}
\begin{split}
&\F(\bar t^{j+1}_{\segg{s_{j+1}}}) \cdot\Pfm\sk{
\F^{j}_{\bar s'}({\bar t}^{j}_{\segg{s_{j}}})}=\prod_{i=1}^{j}\sk{
\prod_{s_{i-1}<\ell<\ell'\leq s_i}
\alpha_3(t^j_{\ell}/t^j_{\ell'})^{-1} }\sum_{A_j}\cdots \sum_{A_2}\sum_{A_1}\\
&\quad \times
Y(t^{j+1}_{a_{s_j}},\ldots,  t^{j+1}_{a_{s_1}};t^j_{s_j},\ldots,t^j_{1})
 \   \prod_{i=1}^j \prod_{\ell<\ell'\atop \ell,\ell'\in A_i}
\!\alpha_1(t^{j+1}_{\ell}/t^{j+1}_{\ell'}) \!\!\!\!\!
\prod_{\ell<\ell'\atop \ell\in \coA_{i-1}\ \ell'\in A_i}
\!\!\! \alpha_2(t^{j+1}_{\ell}/t^{j+1}_{\ell'})\\
&\quad \times
\mathop{\underbrace{F_{j+1}(t^{j+1}_{s_{j+1}})\cdots F_{j+1}(t^{j+1}_{1})}}
\limits_{\
 F_{j+1}(t^{j+1}_{a_1}),\ldots,F_{j+1}(t^{j+1}_{a_{s_j}})\ {\rm are\ omitted}}
\prod^{\longleftarrow}_{1\leq i\leq j}
\sk{ F_{j+2,i}(t^{j+1}_{a_{s_{i-1}+1}})\cdots F_{j+2,i}(t^{j+1}_{a_{s_i}})}.
\end{split}
\end{equation*}
Let us also decompose the summation over the non-ordered sets $A_i=\{a_{s_{i-1}+1},\ldots,a_{s_i}\}$
to the summations over ordered sets $\bar A_i=\{a_{s_{i-1}+1}<\cdots<a_{s_i}\}$ and to the sums over
all permutations among fixed $\{a_{s_{i-1}+1},\ldots,a_{s_i}\}$. Denote the sum over permutation
of the fixed  elements $\{a_{s_{i-1}+1},\ldots,a_{s_i}\}$ as $ \sum_{{\rm per}\bar A_i}$.
The previous formula can be written in the form modulo $\Pfm(U^\coun_{1,j})\cdot U_{1,j+1}$:
\begin{equation}\label{ca13}
\begin{split}
&\F(\bar t^{j+1}_{\segg{s_{j+1}}}) \cdot\Pfm\sk{
\F^{j}_{\bar s'}({\bar t}^{j}_{\segg{s_{j}}})}=\prod_{i=1}^{j}\sk{
\prod_{s_{i-1}<\ell<\ell'\leq s_i}
\alpha_3(t^j_{\ell}/t^j_{\ell'})^{-1}}
\sum_{\bar A_j}\cdots\sum_{\bar A_1} \\
&\!\!\!\!\! \times \prod_{i=1}^j \prod_{\ell<\ell'\atop \ell,\ell'\in \bar A_i}
\!\alpha_1(t^{j+1}_{\ell}/t^{j+1}_{\ell'})
\!\!\!\!\!\!\!  \prod_{\ell<\ell'\atop \ell\in \coA_{i-1}\ \ell'\in \bar A_i}
\!\!\!\!\! \alpha_2(t^{j+1}_{\ell}/t^{j+1}_{\ell'})
\sum_{{\rm per}\bar A_j}\cdots \sum_{{\rm per}\bar A_1}
 Y(t^{j+1}_{a_{s_j}},\ldots,  t^{j+1}_{a_{1}};\bar t^j_{\segg{s_j}})
\\
&\quad \times
\mathop{\underbrace{F_{j+1}(t^{j+1}_{s_{j+1}})\cdots F_{j+1}(t^{j+1}_{1})}}
\limits_{
 F_{j+1}(t^{j+1}_{a_1}),\ldots,F_{j+1}(t^{j+1}_{a_{s_j}})\ {\rm  omitted}}
\prod^{\longleftarrow}_{1\leq i\leq j}
\sk{ F_{j+2,i}(t^{j+1}_{a_{s_{i-1}+1}})\cdots F_{j+2,i}(t^{j+1}_{a_{s_i}})}\ .
\end{split}
\end{equation}

The summations $ \sum_{{\rm per}\bar A_i}$ can be translated to the
$q$-symmetrization of the series $Y(t^{j+1}_{a_{s_j}},\ldots,  t^{j+1}_{a_{1}};t^j_{s_j},\ldots,t^j_{1})$
over the sets of variables $\{t^{j}_{s_{i-1}+1},\ldots,t^{j}_{s_i}\}$ for $i=1,\ldots,j$.
We have the identity of the formal series proved in \cite{KP}:
\begin{align}\notag &
\prod_{i=1}^{j}
\prod_{s_{i-1}<\ell<\ell'\leq s_i}\!\!\!
\alpha_3^{-1}(t^j_{\ell}/t^j_{\ell'})  \prod_{\ell<\ell'\atop \ell,\ell'\in \bar A_i}
\!\alpha_1(t^{j+1}_{\ell}/t^{j+1}_{\ell'})\
\qSym_{\ \bar A_j}\cdots \qSym_{\ \bar A_1}
 Y(t^{j+1}_{a_{s_j}},\ldots,  t^{j+1}_{a_{1}};\bar t^j_{\segg{s_j}})\\
 \label{ca14}
&
\qquad =\ \qSym_{\ \bar t^j_{\seg{s_j}{s_{j-1}}}}\!\!\!\cdots\,
\qSym_{\ \bar t^j_{\seg{s_2}{s_{1}}}}
\qSym_{\ \bar t^j_{\segg{s_1}}}
\sk{ Y(t^{j+1}_{a_{s_j}},\ldots,  t^{j+1}_{a_{1}};\bar t^j_{\segg{s_j}})}.
\end{align}
This identity allows to present  the r.h.s. of \r{ca13} in the form
\begin{equation}\label{ca15}
\begin{split}
&\sum_{\bar A_j}\cdots \sum_{\bar A_2}\sum_{\bar A_1}
\ \     \prod_{\ell<\ell'\atop \ell\in \coA_{i-1}\ \ell'\in \bar A_i}
\!\!\! \frac{q^{-1}-qt^{j+1}_{\ell}/t^{j+1}_{\ell'}}
{q-q^{-1}t^{j+1}_{\ell}/t^{j+1}_{\ell'}}\\
&\quad \times
\mathop{\underbrace{F_{j+1}(t^{j+1}_{s_{j+1}})\cdots F_{j+1}(t^{j+1}_{1})}}
\limits_{
F_{j+1}( t^{j+1}_{a_1}),\ldots,F_{j+1}(t^{j+1}_{a_{s_j}})\ {\rm are\  omitted}}
\prod^{\longleftarrow}_{1\leq i\leq j}
\sk{ F_{j+2,i}(t^{j+1}_{a_{s_i}})\cdots F_{j+2,i}(t^{j+1}_{a_{s_{i-1}+1}})}\\
&\quad\times\ \qSym_{\ \bar t^j_{\seg{s_j}{s_{j-1}}}}\cdots
\qSym_{\ \bar t^j_{\seg{s_2}{s_{1}}}}
\qSym_{\ \bar t^j_{\segg{s_1}}}
\sk{ Y(t^{j+1}_{a_{s_j}},\ldots,  t^{j+1}_{a_{1}};t^j_{s_j},\ldots,t^j_{1})}.
\end{split}
\end{equation}
In its turn, the summation over ordered sets $\bar A_i$ in \r{ca15} can be written as
$q$-symmetri\-za\-tion over the set of variables $\{t^{j+1}_{\segg{s_{j+1}}}\}$, so finally
we obtain the statement of Proposition~\ref{main-fact}.
\hfill$\square$
\medskip

Let $\bar t_{\segg{\bar s}}$ and
$\bar t_{\segg{\bar s'}}$
be the sets of variables
defined by collections of segments $\seg{\bar s}{\bar 0}$ and
$\seg{\bar s'}{\bar 0}$ respectively. Proposition~\ref{main-fact} implies the following
\begin{proposition}\label{fact3}
\begin{equation}\label{pr-ap1}
\begin{split}
&\F(\bar t^{j+1}_{\segg{s_{j+1}}}) \cdot
\qSym_{\ {\bar t}_{\segg{\bar s'}}}
\sk{X(\bar t_{\segg{\bar s'}})\cdot\Pfm\sk{
\F^{j}_{\bar s'}({\bar t}^{j}_{\segg{s_{j}}})}}
=\\
&\qquad=\frac{1}{(s_{j+1}-s_{j})!}\
\qSym_{\ \bar t_{\segg{\bar s}}}
\sk{X(\bar t_{\segg{\bar s}})\cdot
\F^{j+1}_{\bar s}({\bar t}^{j+1}_{\segg{s_{j+1}}})}\quad
{\rm mod}\ \Pfm\sk{U^\coun_{j}}\cdot U_{j+1}
\end{split}
\end{equation}
\end{proposition}

Note that the $q$-symmetrization $\qSym_{\ {\bar t}_{\segg{\bar s'}}}$ in \r{pr-ap1} goes
over the set of variables ${\bar t}_{\segg{\bar s'}}$ which does not include the variables
$\bar t^{j+1}_{\segg{s_{j+1}}}$. This means that the left hand side of \r{pr-ap1} can be written in the
form
$\qSym_{\ {\bar t}_{\segg{\bar s'}}}
\sk{X(\bar t_{\segg{\bar s'}})\cdot\F(\bar t^{j+1}_{\segg{s_{j+1}}}) \cdot\Pfm\sk{
\F^{j}_{\bar s'}({\bar t}^{j}_{\segg{s_{j}}})}}.$
To prove  Proposition~\ref{fact3} we have to substitute \r{ca16} into
this expression. We need the following
\begin{lemma}
Rational series  $\qSym_{\ \bar t^{j-1}_{\segg{s_{j-1}}}}\!\!\cdots
\qSym_{\ \bar t^2_{\segg{s_2}}}
\qSym_{\ \bar t^1_{\segg{s_1}}} X(\bar t_{\bar s'})$ is symmetric in each
group of variables $\{t^j_{s_{i-1}+1},\ldots,t^j_{s_i}\}$, $i=1,\ldots,j$.
\end{lemma}
{\it Proof}\ \ of this Lemma results from the definition of the rational series
\r{rat-Y} and from the following corollary of \rf{ca14}:  the $q$-symmetrization
$\qSym_{\bar v}Y(\bar u;\bar v)$ is a symmetric series on the set
of variables $\bar u$.\hfill$\square$
\smallskip

Statement of Proposition~\ref{fact3} follows now from the identity for the series
\begin{align*}\notag
&\qSym_{\ \bar t^{j}_{\segg{s_{j}}}}\left(
\qSym_{\ \bar t^{j-1}_{\segg{s_{j-1}}}}\cdots
\qSym_{\ \bar t^1_{\segg{s_1}}} X(\bar t_{\bar s'})\times\right.\\
&\left.\quad\times \notag
\qSym_{\ \bar t^j_{\seg{s_j}{s_{j-1}}}}\cdots
\qSym_{\ \bar t^j_{\seg{s_2}{s_{1}}}}
\qSym_{\ \bar t^j_{\segg{s_1}}}
Y(\bar t^{j+1}_{\segg{s_j}}; \bar t^j_{\segg{s_j}})\right)=\\
&\qquad =\prod_{i=1}^j(s_i-s_{i-1})!\ \
\qSym_{\ \bar t_{\segg{\bar s'}}} \sk{X(\bar t_{\bar s})}. 
\end{align*}
\vskip -1cm
\hfill$\square$

\section*{Acknowledgement}
The authors thank S.Loktev, A.Molev and V.Tarasov for useful discussions.
The work of the first author was supported by  RFBR grant
07-02-00878,
 ANR project GIMP No. ANR-05-BLAN-0029-01 and Federal
atomic agency of Russian Federation.
 This work was partially done when the second
author visited Laboratoire d'Annecy-Le-Vieux de Physique Th\'eorique.
He thanks LAPTH for the hospitality and stimulating scientific atmosphere.
His work was supported in part by RFBR grant 06-02-17383 and grant for support
of scientific schools NSh-8065.2006.2. 

\end{document}